\documentclass[12pt,oneside,reqno]{amsart}
\usepackage{amssymb}
\usepackage{}
\usepackage{txfonts}
\usepackage{bbm}
\usepackage{amsmath}
\usepackage{graphicx}
\usepackage{mathrsfs}
\usepackage{stmaryrd}
\usepackage{amsfonts}

\usepackage{enumerate,amsmath,amssymb,amsthm}

\numberwithin{equation}{section}

\newcommand{\be}{\begin{eqnarray}}
\newcommand{\ee}{\end{eqnarray}}
\newcommand{\ce}{\begin{eqnarray*}}
\newcommand{\de}{\end{eqnarray*}}

\newtheorem{theorem}{Theorem}[section]
\newtheorem{lemma}[theorem]{Lemma}
\newtheorem{remark}[theorem]{Remark}
\newtheorem{definition}[theorem]{Definition}
\newtheorem{proposition}[theorem]{Proposition}
\newtheorem{Examples}[theorem]{Example}
\newtheorem{corollary}[theorem]{Corollary}

\def\e{{\mathrm{e}}}
\def\eps{\varepsilon}

\def\p{\partial}

\def\[{{\Big[}}
\def\]{{\Big]}}
\def\<{{\langle}}
\def\>{{\rangle}}
\def\({{\Big(}}
\def\){{\Big)}}

\def\bx{{\mathbf{x}}}

\def\dif{{\mathord{{\rm d}}}}

\def\min{{\mathord{{\rm min}}}}

\def\no{\nonumber}
\def\={&\!\!=\!\!&}
\def\bt{\begin{theorem}}
\def\et{\end{theorem}}
\def\bl{\begin{lemma}}
\def\el{\end{lemma}}
\def\br{\begin{remark}}
\def\er{\end{remark}}

\def\wt{\widetilde}

\def\bd{\begin{definition}}
\def\ed{\end{definition}}
\def\bp{\begin{proposition}}
\def\ep{\end{proposition}}
\def\bc{\begin{corollary}}
\def\ec{\end{corollary}}
\def\bx{\begin{Examples}}
\def\ex{\end{Examples}}

\def\cB{{\mathcal B}}

\def\mE{{\mathbb E}}

\def\mH{{\mathbb H}}

\def\mN{{\mathbb N}}

\def\mP{{\mathbb P}}

\def\mR{{\mathbb R}}

\def\mW{{\mathbb W}}

\def\sF{{\mathscr F}}

\def\sL{{\mathscr L}}

\def\geq{\geqslant}
\def\leq{\leqslant}

\allowdisplaybreaks

\begin{document}

\title{Heat kernels and analyticity of non-symmetric
\\ jump
diffusion semigroups}

\date{}
\author{{Zhen-Qing Chen} \quad
and \quad {Xicheng Zhang}}

\address{Zhen-Qing Chen:
Department of Mathematics, University of Washington, Seattle, WA 98195, USA\\
Email: zqchen@uw.edu
 }
\address{Xicheng Zhang:
School of Mathematics and Statistics, Wuhan University,
Wuhan, Hubei 430072, P.R.China\\
Email: XichengZhang@gmail.com
 }

\begin{abstract}
Let $d\geq 1$ and $\alpha \in (0, 2)$.
Consider the following non-local and non-symmetric L\'evy-type operator
on $\mR^d$:
$$
\sL^\kappa_{\alpha}f(x):=\mbox{p.v.}\int_{\mR^d}(f(x+z)-f(x))
\frac{\kappa(x,z)}{ |z|^{d+\alpha}} \dif z,
$$
where $0<\kappa_0\leq \kappa(x,z)\leq \kappa_1$,
$\kappa(x,z)=\kappa(x,-z)$,  and $|\kappa(x,z)-\kappa(y,z)|\leq\kappa_2|x-y|^\beta$ for some $\beta\in(0,1)$.
Using Levi's method, we construct
the fundamental solution (also called heat kernel)
$p^\kappa_\alpha (t, x, y)$
 of $\sL^\kappa_\alpha $, and
establish its sharp two-sided estimates as well as its
fractional derivative and gradient estimates of the heat kernel.
We also show that $p^\kappa_\alpha (t, x, y)$ is jointly H\"older continuous
in $(t, x)$.
The lower bound heat kernel estimate is obtained by using a probabilistic argument. The fundamental solution of $\sL^\kappa_{\alpha}$ gives rise
a Feller process $\{X, \mP_x, x\in \mR^d\}$ on $\mR^d$. We determine the L\'evy
system of $X$ and show that $\mP_x$ solves the martingale problem
for $(\sL^\kappa_{\alpha}, C^2_b(\mR^d))$.
Furthermore, we  obtain the analyticity of the non-symmetric semigroup associated with $\sL^\kappa_\alpha $ in $L^p$-spaces
for every $p\in[1,\infty)$. A maximum principle for solutions
of the parabolic equation $\partial_t u =\sL^\kappa_\alpha u$ is also established.

\bigskip


\end{abstract}

\thanks{The research of ZC is partially supported
by NSF Grants DMS-0906743 and  NNSFC Grant 11128101.
The research of XZ is partially supported by NNSF grant of China (No. 11271294).}

\maketitle

\section{Introduction}

Let $\sL$ be a second order elliptic differential operator in $\mR^d$ given by
\begin{equation}\label{e:1.1}
\sL f(x)=\sum_{i,j=1}^d \p_i \left( a_{ij}(x)\, \p_jf(x)\right)
 + \sum_{i=1}^d b_i(x) \partial_i  f(x),
\end{equation}
where $(a^{ij}(x))_{1\leq i, j\leq d}$ is a bounded measurable (not necessarily symmetric) $d\times d$-matrix-valued function on $\mR^d$ that is
uniformly elliptic, and $b_i(x)$, $1\leq i\leq d$, are bounded measurable
functions on $\mR^d$. Here $\partial_i f(x)$ stands for the partial derivative $\frac{\partial f(x)}{\partial x_i}$. It is well known that there is a diffusion process $X$
having $\sL$ as its infinitesimal generator; see \cite{CZ}.
The celebrated DeGiorgi-Nash-Moser-Aronson theory asserts
 that every bounded parabolic function of $\sL$ (or
equivalently, of $X$) is locally H\"older continuous and the parabolic Harnack inequality holds for non-negative parabolic functions of $\sL$.
Moreover, $\sL$ has a jointly continuous heat kernel
(or equivalently, transition density function of $X$) $p(t, x, y)$
with respect to the Lebesgue measure on $\mR^d$ that enjoys the
Aronson's Gaussian type estimates.

Quite a lot progress has been made during the last decade
in developing DeGiorgi-Nash-Moser-Aronson type theory for
symmetric non-local operators; see, e.g.,
 \cite{BL2, CaS1, Ch-Ku, Ch-Ku1, Ch}
and the references therein. In particular, it is shown in
Chen and Kumagai \cite{Ch-Ku} that, for every $0<\alpha<2$ and for
any symmetric measurable function
 $c(x, y)$ on $\mR^d\times \mR^d$ that is bounded between two
 positive constants $\kappa_0$ and $\kappa_1$,
 the symmetric non-local operator
 \begin{equation}\label{e:1.2}
 \sL f(x)=\lim_{\eps \to 0} \int_{\{ y\in \mR^d: |y-x|>\eps\}}
 (f(y)-f(x)) \frac{c(x, y)}{|x-y|^{d+\alpha}}
 \dif y
 \end{equation}
 defined in the distributional sense admits a jointly H\"older
 continuous heat kernel $p(t, x, y)$ with respect to the Lebesgue
 measure on $\mR^d$, which satisfies
 \begin{equation}\label{e:1.3}
C^{-1}\,  \frac{t}{(t^{1/\alpha} + |x-y|)^{d+\alpha}}
\leq p(t, x, y)   \leq
C\, \frac{t}{(t^{1/\alpha} + |x-y|)^{d+\alpha}}
\end{equation}
for every $t>0$ and  $x, y\in \mR^d$, where $C\geq 1$
is a constant that depends only on
$d, \alpha, \kappa_0$ and $\kappa_1$.
The operator $\sL$ in \eqref{e:1.2} is symmetric in the sense that
$$
\int_{\mR^d}g(x)\sL  f(x)\dif x=\int_{\mR^d}f(x)\sL  g(x)\dif x
\qquad \hbox{for } f,g\in
C^\infty_c(\mR^d),
$$
where $C^\infty_c(\mR^d)$ denotes the space of smooth functions on $\mR^d$ with
compact support.
When $c(x, y)$ is a positive constant, $\sL$ above
is a constant multiple of the fractional Laplacian
$\Delta^{\alpha/2}:=- (-\Delta)^{\alpha/2}$ on $\mR^d$,
which is the infinitesimal generator of a (rotationally)
symmetric $\alpha$-stable process on $\mR^d$.
The symmetric  non-local stable-like operator $\sL$
defined by \eqref{e:1.2} is the analog to
$\Delta^{\alpha/2}$ of a symmetric uniformly elliptic
divergence
form operator to Laplacian
$\Delta$. Estimate \eqref{e:1.3} can be viewed as
an Aronson type estimate for symmetric stable-like operator
$\sL$ of \eqref{e:1.2}.

The purpose of this paper is to study heat kernels
and their sharp two-sided estimates for non-symmetric
and non-local stable-like operators
of the following form:
\begin{align}
\sL^\kappa_\alpha f(x)
:=\mbox{p.v.}\int_{\mR^d}(f(x+z)-f(x))\frac{\kappa(x,z)}{|z|^{d+\alpha}}
\dif z,\label{Non}
\end{align}
where p.v. stands for the Cauchy principle value; that is
$$
\sL^\kappa_\alpha f(x)=\lim_{\eps \to 0} \int_{\{z\in \mR^d: |z-x|>\eps\}}
(f(x+z)-f(x))\frac{\kappa(x,z)}{|z|^{d+\alpha}}
\dif z.
$$
 Here $d\geq 1$, $0<\alpha<2$,
 and $\kappa(x,z)$ is a measurable function on $\mR^d\times\mR^d$ satisfying
\begin{align}
0<\kappa_0\leq \kappa(x,z)\leq \kappa_1, \qquad
\kappa(x,z)=\kappa(x,-z),    \label{Con1}
\end{align}
and for some $\beta\in(0,1)$
\begin{align}
|\kappa(x,z)-\kappa(y,z)|\leq \kappa_2|x-y|^\beta.\label{Con2}
\end{align}
That $\kappa (x, z)$ is symmetric in $z$
is a commonly assumed condition in the literature of non-local operators;
see \cite{CaS1} for example.
Due to this symmetry condition,   we may write
$$
\sL^\kappa_\alpha f(x)=\frac{1}{2}\int_{\mR^d}(f(x+z)+f(x-z)-2f(x))
\frac{\kappa(x,z)}{|z|^{d+\alpha}}\dif z.
$$
We point out here that, unlike the operator $\sL$ of \eqref{e:1.2},
the operator  $\sL^\kappa_\alpha$ defined by \eqref{Non} is typically
non-symmetric. The relation between $\sL^\kappa_\alpha$ of \eqref{Non} to
$\sL$ of \eqref{e:1.2} is analogous to that of elliptic operators of non-divergence form to elliptic operators  of divergence form.

The following is the main result of this paper.

\bt\label{Main}
Under (\ref{Con1}) and (\ref{Con2}), there exists a unique nonnegative continuous function $p^\kappa_{\alpha} (t,x,y)$ on $(0,1]\times\mR^d\times\mR^d$ solving
\begin{align}
\p_tp^\kappa_{\alpha} (t,x,y)=\sL^\kappa_\alpha p^\kappa_{\alpha}
(t,\cdot,y)(x) ,   \label{eq15}
\end{align}
and satisfying the following four properties:
\begin{enumerate}[(i)]
\item  (Upper bound)
There is a constant $c_1>0$
so that for all $t\in(0,1]$ and $x,y\in\mR^d$,
\begin{align}
p^\kappa_{\alpha} (t,x,y)\leq c_1 t ( t^{1/\alpha}+|x-y|)^{-d-\alpha}. \label{eq16}
\end{align}
\item (H\"older's estimate)
For every $\gamma \in(0,\alpha\wedge 1)$, there is a constant
$c_2>0$ so that
for every
$t\in(0,1]$ and $x,x',y\in\mR^d$,
\begin{align}
|p^\kappa_{\alpha} (t,x,y)-p^\kappa_{\alpha} (t,x',y)|\leq c_2|x-x'|^\gamma t^{1-\frac{\gamma}{\alpha}}\left(t^{1/\alpha}+|x-y|\wedge|x'-y|)^{-d-\alpha} \right). \label{Ho}
\end{align}
\item (Fractional derivative estimate) For all $x,y\in\mR^d$, the mapping $t\mapsto \sL^\kappa_\alpha p^\kappa_{\alpha}(t,\cdot,y)(x)$ is continuous on $(0,1]$, and
\begin{align}
|\sL^\kappa_\alpha p^\kappa_{\alpha} (t,\cdot,y)(x)|\leq c_3(t^{1/\alpha}+|x-y|)^{-d-\alpha}.\label{eq17}
\end{align}
\item (Continuity) For any bounded and uniformly continuous function $f:\mR^d\to\mR$,
\begin{align}\label{eq14}
\lim_{t\downarrow 0}\sup_{x\in\mR^d}\left|\int_{\mR^d}p^\kappa_{\alpha} (t,x,y)f(y)\dif y-f(x)\right|=0.
\end{align}
\end{enumerate}
Moreover,  we have the following conclusions.
\begin{enumerate}[(1)]
\item The constants $c_1$, $c_2$ and $c_3$ in (i)-(iii) above
can be chosen so
that they  depend only
on $(d, \alpha, \beta, \kappa_0, \kappa_1, \kappa_2)$,  $(d, \alpha, \beta, \gamma, \kappa_0, \kappa_1, \kappa_2)$, and $(d, \alpha, \beta, \kappa_0, \kappa_1, \kappa_2)$, respectively.

\item (Conservativeness)
For all $(t,x)\in(0,1]\times\mR^d$,
$p^\kappa_\alpha (t, x, y)\geq 0$ and
\begin{align}
\int_{\mR^d}p^\kappa_{\alpha} (t,x,y)\dif y=1.\label{EK1}
\end{align}
\item (C-K equation) For all $s,t\in(0,1]$ and $x,y\in\mR^d$, the following Chapman-Kolmogorov's equation holds:
\begin{align}\label{eq21}
\int_{\mR^d}p^\kappa_{\alpha} (t,x,z)p^\kappa_{\alpha} (s,z,y)\dif z=p^\kappa_{\alpha} (t+s,x,y).
\end{align}
\item  (Lower bound) For all $t\in(0,1]$ and $x,y\in\mR^d$,
\begin{align}
p^\kappa_{\alpha} (t,x,y)\geq c_4t(t^{1/\alpha}+|x-y|)^{-d-\alpha}. \label{e:1.14}
\end{align}
\item (Gradient estimate) If $\alpha\in[1,2)$, for all $x,y\in\mR^d$ and $t\in(0,1]$,
\begin{align}
|\nabla p^\kappa_{\alpha} (t,\cdot,y)(x)|\leq c_5t^{1-1/\alpha}(t^{1/\alpha}+|x-y|)^{-d-\alpha}.\label{eq170}
\end{align}
\item (Generator) For any $f\in C^2_b(\mR^d)$, we have
\begin{align}\label{eq233}
\lim_{t\downarrow 0}\frac{1}{t}\int_{\mR^d}p^\kappa_{\alpha} (t,x,y)(f(y)-f(x))\dif y=\sL^\kappa_\alpha f(x),
\end{align}
and the convergence is uniform.

\item (Analyticity) The $C_0$-semigroup $( P^\kappa_t)_{t\geq 0}$ of $\sL^\kappa_\alpha$ defined
by $P^\kappa_t f(x):=\int_{\mR^d} p^\kappa_\alpha (t, x, y) f(y) \dif y$
is analytic in $L^p(\mR^d)$ for every  $p\in[1,\infty)$.
\end{enumerate}
\et

Here $C^2_b(\mR^d)$ is the space of bounded continuous functions on
$\mR^d$ that have bounded continuous first and second order partial derivatives.
 A $C_0$-semigroup means a strongly continuous semigroup
in the space of continuous functions on $\mR^d$ that vanish
at infinity equipped with the uniform topology.

\br
\eqref{eq16} and \eqref{e:1.14} give the sharp two-sided estimates
of the heat kernel $p^\kappa_\alpha (t, x, y)$.
We can restate Estimate \eqref{eq17} as
$$
 |\partial_t  p^\kappa_{\alpha} (t,x,y) |\leq c_3(t^{1/\alpha}+|x-y|)^{-d-\alpha}.
$$
This together with \eqref{Ho} and (1) of Theorem \ref{Main} yields
that for $0<s<t$ and $x, x', y\in \mR^d$,
\begin{align}
&|p^\kappa_{\alpha} (s,x,y)-p^\kappa_{\alpha} (t,x',y)| \nonumber \\
&\quad\leq \wt c_2 \left( |t-s|+|x-x'|^\gamma t^{1-\frac{\gamma}{\alpha}} \right) \left(s^{1/\alpha}+|x-y|\wedge |x'-y|\right)^{-d-\alpha},
\label{e:1.17}
\end{align}
where
$\wt c_2=c_2+ c_3$.
\er

\medskip

To the best of the authors' knowledge, Theorem \ref{Main} is the first result on heat kernels and their estimates for a general class of non-symmetric and
  non-local stable-like operators under H\"oder continuous condition in
 $x\mapsto \kappa (x, z)$.
 We mention that in the framework of pseudodifferential operator theory, Kochubei \cite{Ko} (see also \cite{Ei-Iv-Ko})
has already studied the existence of fundamental solutions for $\sL^\kappa_\alpha$ by using Levi's method.
But strong smoothness of $\kappa(x,y)$ in $y$ and $\alpha\in[1,2)$ are required.
In Chen and Wang \cite{CW}, fractional Laplacian
$\Delta^{\alpha/2}$ perturbed by lower order non-local operator
is studied, which corresponds to the case when
$\kappa (x, z)= a + b (x, z) |z|^{\alpha -\delta}$
for some constant $a>0$ and a
bounded measurable $b(x, z)$ with  $b(x, z)=b(x, -z)$.
As a special case of the much more general results obtained in \cite{CW},
it is proved there that for this type of $\kappa (x, z)$,
 when there are two positive constants $\kappa_0, \kappa_1$
 so that $\kappa_0 \leq \kappa (x, z)\leq \kappa_1$ (but no H\"older
continuity is assumed in $x\mapsto b(x, z)$),
 $\sL^\kappa_\alpha$ has a unique jointly continuous heat kernel
$p^\kappa_\alpha (t, x, y)$ and it enjoys the two-sided estimates
\eqref{eq16} and \eqref{e:1.14}.

Although quite a lot is known for symmetric non-local operators,
there are very limited results in literature
on heat kernel estimates for non-symmetric and non-local operators. In \cite{Bo-Ja0}, Bogdan and Jakubowski studied the estimates of heat kernel of $\Delta^{\alpha/2}$ perturbed by a gradient operator with $\alpha\in(1,2)$ (see also \cite{Wa-Zh} for some extension).
Jakubowski and Szczypkowski \cite{Ja-Sz2} considered the time-dependent gradient perturbation of $\Delta^{\frac{\alpha}{2}}$, while
 Jakubowski \cite{Ja} established the global time estimate of heat kernel of $\Delta^{\alpha/2}$ under small singular drifts.
In \cite{Ch-Ki-So0, Ch-Ki-So, Ch-Ki-So1}, Chen, Kim and Song obtained sharp two-sided estimates for the Dirichlet heat kernel
of $\Delta^{\alpha/2}$ as well as of its gradient
and Feynman-Kac perturbations.
Global as well as Dirichlet heat kernel estimates for non-local operators
 $\Delta +\Delta^{\alpha/2} + b\cdot \nabla$ and for $m-(m^{2/\alpha}-\Delta)^{\alpha/2}+b\cdot \nabla$
 have been investigated in
 \cite{CH} and \cite{CWa}, respectively.
In the critical case of $\alpha=1$, the sharp two-sided heat kernel estimates of $\Delta^{1/2}+b\cdot \nabla $ with H\"older continuous drift $b$
 was obtained recently in \cite{Xi-Zh}
by using a Levi's method. In \cite{Ma-Mi}, Maekawa and Miura obtained the upper bounds estimates
for the fundamental solutions of general non-local diffusions with divergence free drift.

\medskip

We next briefly describe the approach of this paper.
For the construction and upper bound estimates of the heat kernel,
we use a method based on Levi's freezing coefficients argument (cf. \cite{La-So-Ur, Fr}).
However, in contrast
to the previous work \cite{Xi-Zh},
 a new way to freeze the
coefficient $\kappa (x, z)$ is needed (see Section 3).
This causes quite many new challenges.
In particular, we need to estimate the fractional derivative
of the freezing heat kernel
and to prove the continuous dependence of heat kernels with respect to the kernel function $\kappa$ (see Subsections 2.3 and 2.4).
Strong stability of the heat kernels in terms of
the maximal distance between jumping kernels has recently been studied
in Bass and Ren \cite{BR} (see Theorem 5.3 there)
for symmetric stable-like operators
\eqref{e:1.2}. But here we need a more refined stability results on
the heat kernels and their derivatives; see Theorem \ref{T:2.5} below.
To show the uniqueness and non-negativeness of the heat kernel,
we establish a maximum principle
for solutions of the parabolic equation
$\partial_t u(t, x)=\sL^\kappa_\alpha u(t, x)$; see Theorem \ref{Th1}.
For the lower bound estimate \eqref{e:1.14} on the heat kernel,
we use a probabilistic approach. The heat kernel $p^\kappa_\alpha (t, x, y)$
determines a strong Feller process $X=\{X_t, t\geq 0; \mP_x, x\in \mR^d\}$
on $\mR^d$. We show that for each $x\in \mR^d$, $\mP_x$ solves the martingale
problem for $(\sL^\kappa_\alpha, C^2_b (\mR^d))$ with initial value $x$;
see \eqref{ERY1} below.
We then deduce from it the L\'evy system of $X$, which tells us that
$k(x, z)|z|^{-(d+\alpha)}$ is the jump intensity of $X$ making a jump from
$x$ with size $z$. The lower bound estimate for $p^\kappa_\alpha$ can then
be obtained by   a probabilistic argument involving the use of the
L\'evy system of $X$.

\br It
will be shown in a subsequent paper \cite{CZh} that
solution to the martingale problem for $(\sL^\kappa_\alpha, C^\infty_c(\mR^d))$ is unique. (In fact it will be established for a more general class of non-local operators.)
Thus the heat kernel $p^\kappa_\alpha (t, x, y)$ in Theorem \ref{Main}
can also be regarded
as the (unique) transition density function of the unique solution
to the martingale problem for $(\sL^\kappa_\alpha, C^\infty_c(\mR^d))$.
\er

Notion of analyticity of a $C_0$-semigroup plays a central role in the semigroup theory of evolution equations (cf. \cite{Fr0, He, Pa}).
For differential operators $\sL$ of \eqref{e:1.1},
it is well-known that its associated $C_0$-semigroup
is analytic in $L^p$-spaces for every $p\in(1,\infty)$ at least
when $a_{ij}$ are smooth (cf. \cite[Chapter 7]{Pa}).
The proof of this fact is based upon the following deep a priori estimate:
$$
\|\p_i\p_jf\|_p\leq C(\|\sL^a_2  f\|_p+\|f\|_p),\quad f\in\mW^{2,p}(\mR^d),
$$
which is a consequence of singular integral operator theory.
For nonlocal operator $\sL^\kappa_\alpha$ of \eqref{e:1.2},
under some additional assumptions on $\kappa(x,z)$, it was shown in \cite{Zh1} and \cite{Zh2} that for any $p\in(1,\infty)$ and $\alpha\in(0,2)$,
$$
c_6\|f\|_{\mH^{\alpha,p}}\leq \|\sL^\kappa_\alpha f\|_p+\|f\|_p\leq c_6^{-1}\|f\|_{\mH^{\alpha,p}}, \ \ f\in\mH^{\alpha,p},
$$
where $\mH^{\alpha,p}=(I-\Delta)^{-\frac{\alpha}{2}}(L^p)$ is the usual Bessel potential space. In this case, it is possible to show the analyticity of its associated semigroup $(P^\kappa_t)_{t\geq 0}$
by using Agmon's method \cite{Fr0}. However in this paper we are able to establish the analyticity of the semigroup $(P^\kappa_t)_{t\geq 0}$
without these additional assumptions. We achieve this
by establishing the inequality $\| \sL^\kappa_\alpha P_t f\|_p \leq c t^{-1}
\|f\|_p$ for every $t>0$ and $p\in [1, \infty)$.

\medskip

We now give an application of Theorem \ref{Main}
to stochastic differential equations driven by
(rotationally) symmetric stable processes.
Suppose that $A(x)=(a_{ij}(x))_{1\leq i, j\leq d}$
is a bounded continuous $d\times d$-matrix-valued function
on $\mR^d$ that is nondegenrate at every $x\in \mR^d$,
and $Y$ is a (rotationally) symmetric $\alpha$-stable process
on $\mR^d$ for some $0<\alpha <2$. It is shown in Bass and Chen
\cite[Theorem 7.1]{BC} that for every $x\in \mR^d$, SDE
\begin{equation}\label{e:1.18}
\dif X_t = A(X_{t-}) \dif Y_t, \qquad X_0=x,
\end{equation}
has a unique weak solution.
(Although in \cite{BC} it is assumed $d\geq 2$, the
 results there are valid for $d=1$ as well.)
The family of these weak solutions
forms a strong Markov process
$\{X, \mP_x, x\in \mR^d\}$.
Using It\^o's
formula, one deduces (see the display above (7.2) in
\cite{BC}) that $X$ has   generator
\begin{equation}\label{e:1.19}
\sL f(x) = {\rm p.v.} \int_{\mR^d} \left( f(x+A(x)y)-f(x)\right)
\frac{c_{d, \alpha}} {|y|^{d+\alpha}} \dif y,
\end{equation}
where $c_{d, \alpha}$ is a positive constant that depends on $d$ and $\alpha$.
A change of variable formula $z=A(x) y$ yields
\begin{equation}\label{e:1.20}
\sL f(x) = {\rm p.v.} \int_{\mR^d} \left( f(x+z)-f(x)\right)
\frac{\kappa (x, z)} {|z|^{d+\alpha}} \dif z ,
\end{equation}
where
\begin{equation}\label{e:1.21}
\kappa (x, z)= \frac{c_{d,\alpha}}{|{\rm det} A(x)|} \left( \frac{|z|}{|A(x)^{-1}z|}\right)^{d+\alpha}.
\end{equation}
Here ${\rm det}(A(x))$ is the determinant of the matrix $A(x)$
and $A(x)^{-1}$ is the inverse of $A(x)$.
As an application of the main result of this paper, we have

\bc \label{C:1.4} Suppose that $A(x)=(a_{ij}(x))$ is uniformly bounded and elliptic (that is, there are positive constants
$\lambda_0$ and $\lambda_1$ so that
 $\lambda_0 I_{d\times d} \leq A(x)\leq \lambda_1 I_{d\times d}$ for every
 $x\in \mR^d$) and there are $\beta\in (0, 1)$ and $\lambda_2>0$
 so that
 $$ |a_{ij}(x)-a_{ij}(y)| \leq \lambda_2 |x-y|^\beta
 \quad \hbox{for } 1\leq i, j\leq d.
 $$
Then the strong Markov process $X$ formed by the unique weak solution to SDE \eqref{e:1.18} has a jointly continuous transition density function
$p(t, x, y)$ with respect to the Lebesgue measure on $\mR^d$ and
there is a constant $C>0$ that depends only on $(d, \alpha, \beta, \lambda_0, \lambda_1)$ so that
$$
C^{-1}\,  \frac{t}{(t^{1/\alpha} + |x-y|)^{d+\alpha}}
\leq p(t, x, y)   \leq
C\, \frac{t}{(t^{1/\alpha} + |x-y|)^{d+\alpha}}
$$
for every $t\in (0, 1]$ and $x, y\in \mR^d$.
\ec

\medskip

The remainder of the paper is organized as follows.
 In Section 2, we prepare some necessary results about the estimates of the heat kernel of {\it spatial-independent} symmetric L\'evy operators.
In Section 3, we construct the heat kernel of {\it spatial-dependent} L\'evy operators by using
Levi's method. Lastly, in Section 4 we present
the proof of  the main result of this paper, Theorem \ref{Main}.

We conclude this section by introducing the following conventions.
 The letter $C$ with or without subscripts
will denote a positive constant, whose value is not important and may change in different places. We write $f(x)\preceq g(x)$ to mean that there exists a constant
$C_0>0$ such that $f(x)\leq C_0 g(x)$; and $f(x)\asymp g(x)$ to mean
that there exist $C_1,C_2>0$ such that $C_1 g(x)\leq f(x)\leq C_2 g(x)$.
We will also use the abbreviation $f(x\pm z)$ for $f(x+z)+f(x-z)$.
For $p\geq 1$, $L^p$-norm of $L^p(\mR^d)=L^p(\mR^d; \dif x)$ will be denoted as $\| f\|_p$.
We use ``$:=$'' to denote a definition. For $a, b\in \mR$, $a\wedge b:= \min\{a, b\}$ and $a\vee b:=\max\{a, b\}$.

\section{Preliminaries}

Throughout this paper, we shall fix $\alpha\in(0,2)$ and assume
$$
(t,x)\in(0,1]\times\mR^d.
$$
For $\gamma,\beta\in\mR$, we introduce the following function on $(0,1]\times\mR^d$ for later use:
\begin{align}
\varrho^\beta_\gamma(t,x):=t^{\frac{\gamma}{\alpha}}(|x|^\beta\wedge 1)(t^{1/\alpha}+|x|)^{-d-\alpha}.\label{Def2}
\end{align}
\subsection{Convolution inequalities}
The following lemma will play an important role in the sequel, which is similar to \cite[Lemma 1.4]{Kol}
and \cite[Lemma 2.3]{Xi-Zh}.

\bl\label{Le09}
(i) For all $\beta\in[0,\frac{\alpha}{2}]$ and $\gamma\in\mR$, we have
\begin{align}
\int_{\mR^d}\varrho^\beta_\gamma(t,x)\dif x\preceq  t^{\frac{\gamma+\beta-\alpha}{\alpha}},\ \ (t,x)\in(0,1)\times\mR^d.\label{ES4}
\end{align}
(ii) For all $\beta_1,\beta_2\in[0,\frac{\alpha}{4}]$ and $\gamma_1,\gamma_2\in\mR$, we have
\begin{align}
\int_{\mR^d}\varrho^{\beta_1}_{\gamma_1}(t-s,x-z)\varrho^{\beta_2}_{\gamma_2}(s,z)
\dif z
\preceq & \left( (t-s)^{\frac{\gamma_1+\beta_1+\beta_2-\alpha}{\alpha}}
s^{\frac{\gamma_2}{\alpha}}+(t-s)^{\frac{\gamma_1}{\alpha}}
s^{\frac{\gamma_2+\beta_1+\beta_2-\alpha}{\alpha}}\right) \varrho^0_0(t,x)\no\\
&+(t-s)^{\frac{\gamma_1+\beta_1-\alpha}{\alpha}}s^{\frac{\gamma_2}{\alpha}}\varrho^{\beta_2}_0(t,x)+(t-s)^{\frac{\gamma_1}{\alpha}}s^{\frac{\gamma_2+\beta_2-\alpha}{\alpha}}\varrho^{\beta_1}_0(t,x).\label{EU7}
\end{align}
(iii) If $\gamma_1+\beta_1>0$ and $\gamma_2+\beta_2>0$, then
\begin{align}
&\int^t_0\!\!\!\int_{\mR^d}\varrho^{\beta_1}_{\gamma_1}(t-s,x-z)\varrho^{\beta_2}_{\gamma_2}(s,z)\dif z\dif s\no\\
&\qquad \preceq \cB(\tfrac{\gamma_1+\beta_1}{\alpha},\tfrac{\gamma_2+\beta_2}{\alpha})
\left( \varrho^0_{\gamma_1+\gamma_2+\beta_1+\beta_2}
+\varrho^{\beta_1}_{\gamma_1+\gamma_2+\beta_2}
+\varrho^{\beta_2}_{\gamma_1+\gamma_2+\beta_1}\right)(t,x), \label{eq30}
\end{align}
where $\cB(\gamma,\beta)$ is the usual Beta function defined by
$$
\cB(\gamma,\beta):=\int^1_0(1-s)^{\gamma-1}s^{\beta-1}\dif s,\ \ \gamma,\beta>0.
$$
Moreover, the constants contained in the above $\preceq$ only depend on
$d$, $\alpha$, and the $\beta$'s.
\el

\begin{proof}
(i) Notice that
\begin{align*}
\int_{\mR^d}\frac{|x|^\beta}{(t^{1/\alpha}+|x|)^{d+\alpha}}\dif x
&\preceq\int^\infty_0\frac{r^{\beta+d-1}}{(t^{1/\alpha}+r)^{d+\alpha}}\dif r
=\left(\int^{t^{1/\alpha}}_0+\int^\infty_{t^{1/\alpha}}\right)\frac{r^{\beta+d-1}}{(t^{1/\alpha}+r)^{d+\alpha}}\dif r\\
&\leq \int^{t^{1/\alpha}}_0\frac{r^{d+\beta-1}}{t^{(d+\alpha)/\alpha}}\dif r+\int^\infty_{t^{1/\alpha}}r^{\beta-1-\alpha}\dif r
=\frac{t^{(\beta-\alpha)/\alpha}}{d+\beta}+\frac{t^{(\beta-\alpha)/\alpha}}{\alpha-\beta},
\end{align*}
which implies (\ref{ES4}) by definition.

(ii) In view of
$$
(t^{1/\alpha}+|x|)^{d+\alpha}\leq C_{d,\alpha}\left(((t-s)^{1/\alpha}+|x-z|)^{d+\alpha}+(s^{1/\alpha}+|z|)^{d+\alpha}\right),
$$
we have
\begin{align}
\varrho^0_0(t-s,x-z)\varrho^0_0(s,z)\leq
C_{d,\alpha}\left(\varrho^0_0(t-s,x-z)+\varrho^0_0(s,z)\right)
\varrho^0_0(t,x).\label{ET2}
\end{align}
Noticing that by $(a+b)^\beta\leq a^\beta+b^\beta$ for $\beta\in(0,1)$,
\begin{align*}
(|x-z|^{\beta_1}\wedge 1)(|z|^{\beta_2}\wedge 1)
&\leq(|x-z|^{\beta_1}\wedge 1)((|x-z|^{\beta_2}+|x|^{\beta_2})\wedge 1)\\
&\leq|x-z|^{\beta_1+\beta_2}\wedge 1+(|x-z|^{\beta_1}\wedge 1)(|x|^{\beta_2}\wedge 1),\\
(|x-z|^{\beta_1}\wedge 1)(|z|^{\beta_2}\wedge 1)
&\leq((|z|^{\beta_1}+|x|^{\beta_1})\wedge 1)(|z|^{\beta_2}\wedge 1)\\
&\leq|z|^{\beta_1+\beta_2}\wedge 1+(|x|^{\beta_1}\wedge1)(|z|^{\beta_2}\wedge 1),
\end{align*}
we have
\begin{align*}
&\varrho^{\beta_1}_{\gamma_1}(t-s,x-z)\varrho^{\beta_2}_{\gamma_2}(s,z)
=(t-s)^{\frac{\gamma_1}{\alpha}}s^{\frac{\gamma_2}{\alpha}}(|x-z|^{\beta_1}\wedge 1)(|z|^{\beta_2}\wedge 1)
\varrho^0_0(t-s,x-z)\varrho^0_0(s,z)\no\\
&\quad\preceq(t-s)^{\frac{\gamma_1}{\alpha}}s^{\frac{\gamma_2}{\alpha}}\Big\{|x-z|^{\beta_1+\beta_2}\wedge 1+(|x-z|^{\beta_1}\wedge 1)(|x|^{\beta_2}\wedge 1)\Big\}\varrho^0_0(t-s,x-z)\varrho^0_0(t,x)\no\\
&\qquad+(t-s)^{\frac{\gamma_1}{\alpha}}s^{\frac{\gamma_2}{\alpha}}\Big\{|z|^{\beta_1+\beta_2}\wedge 1+(|x|^{\beta_1}\wedge 1)(|z|^{\beta_2}\wedge 1)\Big\}\varrho^0_0(s,z)\varrho^0_0(t,x)\no\\
&\quad\preceq s^{\frac{\gamma_2}{\alpha}}\Big\{\varrho^{\beta_1+\beta_2}_{\gamma_1}(t-s,x-z)\varrho^0_0(t,x)
+\varrho^{\beta_1}_{\gamma_1}(t-s,x-z)\varrho^{\beta_2}_0(t,x)\Big\}\no\\
&\qquad+(t-s)^{\frac{\gamma_1}{\alpha}}\Big\{\varrho^{\beta_1+\beta_2}_{\gamma_2}(s,z)\varrho^0_0(t,x)
+\varrho^{\beta_2}_{\gamma_2}(s,z)\varrho^{\beta_1}_0(t,x)\Big\}.
\end{align*}
Integrating both sides with respect to $z$ and using (i), we obtain (ii).

(iii) Observe that for $\gamma,\beta>0$,
\begin{align}
\int^t_0(t-s)^{\gamma-1}s^{\beta-1}\dif s=t^{\gamma+\beta-1}\cB(\gamma,\beta).\label{eq32}
\end{align}
Integrating both sides of (\ref{EU7}) with respect to $s$ from $0$ to $t$, we obtain
\begin{align*}
&\int^t_0\!\!\!\int_{\mR^d}\varrho^{\beta_1}_{\gamma_1}(t-s,x-z)\varrho^{\beta_2}_{\gamma_2}(s,z)\dif z\dif s\no\\
&\qquad\preceq t^{\frac{\gamma_1+\gamma_2+\beta_1+\beta_2}{\alpha}}\Big\{\cB(\tfrac{\gamma_1+\beta_1+\beta_2}{\alpha},\tfrac{\gamma_2+\alpha}{\alpha})
+\cB(\tfrac{\gamma_2+\beta_1+\beta_2}{\alpha},\tfrac{\gamma_1+\alpha}{\alpha})\Big\}\varrho^0_0(t,x)\\
&\qquad+t^{\frac{\gamma_1+\gamma_2+\beta_1}{\alpha}}\cB(\tfrac{\gamma_1+\beta_1}{\alpha},\tfrac{\gamma_2+\alpha}{\alpha})\varrho^{\beta_2}_0(t,x)
+t^{\frac{\gamma_1+\gamma_2+\beta_2}{\alpha}}\cB(\tfrac{\gamma_2+\beta_2}{\alpha},\tfrac{\gamma_1+\alpha}{\alpha})\varrho^{\beta_1}_0(t,x),
\end{align*}
which implies  (\ref{eq30}) by $\beta_1,\beta_2<\alpha$ and that $\cB(\gamma,\beta)$ is symmetric and non-increasing with respect to
variables $\gamma$ and $\beta$.
\end{proof}

\subsection{Some estimates of heat kernel of $\Delta^{\frac{\alpha}{2}}$}

Let
$(Z^{(\alpha)}_t)_{t\geq 0}$
be a rotationally invariant $d$-dimensional $\alpha$-stable process, and $p_\alpha(t,x)$
 its probability transition density function with respect to the Lebesgue
 measure on $\mR^d$.
 By the scaling property of $Z^{(\alpha)}_t\stackrel{(d)}{=}t^{1/\alpha}Z^{(\alpha)}_1$,
it is easy to see that
\begin{align}
p_\alpha(t,x)=t^{-d/\alpha}p_\alpha(1,t^{-1/\alpha} x).\label{ER1}
\end{align}
Let $(W_t)_{t\geq 0}$ be a $d$-dimensional standard Brownian motion, and $S^{(\alpha)}_t$ an $\alpha/2$-stable subordinator. It is well-known that $Z^{(\alpha)}_t$ can be realized as
$$
Z^{(\alpha)}_t=W_{S^{(\alpha)}_t}.
$$
Let $\eta_t(s)$ be the density of $S^{(\alpha)}_t$. By subordination, we have
$$
p_\alpha(t,x)=\int^\infty_0(2\pi s)^{-\frac{d}{2}}\e^{-\frac{|x|^2}{2s}}\eta_t(s)\dif s.
$$
By \cite[Theorem 2.1]{Bl-Ge}, one knows that
\begin{align}
p_\alpha(t,x)\asymp \varrho^0_\alpha(t,x)=t(t^{1/\alpha}+|x|)^{-d-\alpha}.\label{ER60}
\end{align}
The following obvious inequality will be used frequently:
\begin{align}
(t^{1/\alpha}+|x+z|)^{-\gamma}\leq 2^\gamma(t^{1/\alpha}+|x|)^{-\gamma},\ \  |z|\leq t^{1/\alpha}\vee(|x|/2).\label{ER5}
\end{align}
Below, for a function $f$ defined on $\mR_+\times\mR^d$, we shall simply write
\begin{equation}\label{e:2.10}
\delta_f(t,x;z):=f(t,x+z)+f(t,x-z)-2f(t,x).
\end{equation}

We need the following lemma.

\bl\label{Le4}
There is a constant $C=C(d, \alpha)>0$ so that for every $t>0$,
$x, x', z\in \mR^d$,
\begin{align}
|\nabla^kp_\alpha(t,x)|
 \leq C \,
t(t^{1/\alpha}+|x|)^{-d-\alpha-k},\ \  k\in\mN,\label{ER00}
\end{align}
\begin{align}
|p_\alpha(t,x)-p_\alpha(t,x')|
\leq C\,
((t^{-1/\alpha}|x-x'|)\wedge 1)\left( p_\alpha(t,x)+p_\alpha(t,x')\right),\label{ER20}
\end{align}
\begin{align}
|\delta_{p_\alpha}(t,x;z)|
\leq C\,
((t^{-\frac{2}{\alpha}}|z|^2)\wedge 1)\left( p_\alpha(t,x\pm z)+p_\alpha(t,x)\right),   \label{ER2}
\end{align}
\begin{align}
&|\delta_{p_\alpha}(t,x;z)-\delta_{p_\alpha}(t,x';z)|
\leq C\,
((t^{-1/\alpha}|x-x'|)\wedge 1)((t^{-\frac{2}{\alpha}}|z|^2)\wedge 1)\no\\
&\qquad\qquad\times\left( p_\alpha(t,x\pm z)+p_\alpha(t,x)+p_\alpha(t,x'\pm z)+p_\alpha(t,x')\right) .  \label{ER3}
\end{align}
\el

\begin{proof}
By the scaling property (\ref{ER1}), it suffices to prove these estimates for $t=1$.

(i) Noticing that (cf. \cite[Theorem 37.1]{De})
$$
\eta_1(s)\preceq s^{-1-\frac{\alpha}{2}}\e^{-s^{-\alpha/2}}\leq s^{-1-\frac{\alpha}{2}},
$$
we have for $|x|>1$,
\begin{align*}
|\nabla p_\alpha(1,x)|\preceq |x|\int^\infty_0s^{-\frac{d}{2}-2-\frac{\alpha}{2}}\e^{-\frac{|x|^2}{2s}}\dif s=|x|^{-d-\alpha-1}\int^\infty_0u^{\frac{d+\alpha}{2}}\e^{-\frac{u}{2}}\dif u.
\end{align*}
Hence,
$$|\nabla p_\alpha(1,x)|\preceq(1+|x|)^{-d-\alpha-1},\ \ x\in\mR^d,
$$
which gives (\ref{ER00}) for $k=1$. The estimates of higher order derivatives are similar.

(ii) Observe that
\begin{align}
p_\alpha(1,x)-p_\alpha(1,x')=\int^1_0\nabla_{x-x'}p_\alpha(1,x+\theta(x'-x))\dif\theta.\label{Tay}
\end{align}
If $|x-x'|\leq 1$, then by (\ref{ER00}), we have
\begin{align*}
| p_\alpha(1,x)-p_\alpha(1,x')|&\preceq |x-x'|\int^1_0(1+|x+\theta(x'-x)|)^{-d-\alpha-1}\dif\theta\\
&\stackrel{(\ref{ER5})}{\preceq} |x-x'|(1+|x|)^{-d-\alpha-1}\stackrel{(\ref{ER60})}{\preceq} |x-x'|p_\alpha(1,x).
\end{align*}
So,
$$
|p_\alpha(1,x)-p_\alpha(1,x')|\preceq (|x-x'|\wedge 1)\Big\{p_\alpha(1,x)+p_\alpha(1,x')\Big\}.
$$
Estimate (\ref{ER20}) follows.

(iii) By using (\ref{Tay}) twice, we have
\begin{align}
\delta_{p_\alpha}(1,x;z)&=p_\alpha(1,x+z)+p_\alpha(1,x-z)-2p_\alpha(1,x)\no\\
&=\int^1_0(\nabla_z p_\alpha(1,x+\theta z)-\nabla_z p_\alpha(1,x-\theta z))\dif \theta\no\\
&=\int^1_0\!\!\!\int^1_0\theta\nabla_z\nabla_z p_\alpha(1,x+(1-2\theta')\theta z)\dif\theta'\dif \theta.\label{ER770}
\end{align}
If $|z|>1$, then
$$
|\delta_{p_\alpha}(1,x;z)|\leq p_\alpha(1,x+z)+p_\alpha(1,x-z)+2p_\alpha(1,x).
$$
If $|z|\leq 1$, then by (\ref{ER00}), we have
\begin{align*}
|\delta_{p_\alpha}(1,x;z)|&\leq|z|^2\int^1_0\!\!\!\int^1_0|\nabla^2 p_\alpha(1,x+(1-2\theta')\theta z)|\dif\theta'\dif \theta\\
&\preceq|z|^2\int^1_0\!\!\!\int^1_0(1+|x+(1-2\theta')\theta z )|)^{-d-\alpha-2}\dif\theta'\dif \theta\\
&\stackrel{(\ref{ER5})}{\preceq}|z|^2(1+|x|)^{-d-\alpha-2}\stackrel{(\ref{ER60})}{\preceq}|z|^2p_\alpha(1,x).
\end{align*}
Hence,
\begin{align}
|\delta_{p_\alpha}(1,x;z)|\preceq(|z|^2\wedge 1)\Big\{p_\alpha(1,x\pm z)+p_\alpha(1,x)\Big\},\label{ER9}
\end{align}
which yields (\ref{ER2}).

(iv) If $|z|\leq 1$ and $|x-x'|\leq 1$, then by (\ref{ER770}) and (\ref{ER00}), we have
\begin{align}
&|\delta_{p_\alpha}(1,x;z)-\delta_{p_\alpha}(1,x';z)|\no\\
&\quad\preceq |x-x'|\cdot|z|^2\int^1_0\!\!\!\int^1_0\!\!\!\int^1_0|\nabla^3p_\alpha|(1,x+(1-2\theta')\theta z+\theta''(x'-x))\dif\theta''\dif\theta'\dif\theta\no\\
&\quad\preceq |x-x'|\cdot|z|^2\int^1_0\!\!\!\int^1_0\!\!\!\int^1_0(1+|x+(1-2\theta')\theta z+\theta''(x'-x)|)^{-d-\alpha-3}\dif\theta''\dif\theta'\dif\theta\no\\
&\quad\stackrel{(\ref{ER5})}{\preceq} |x-x'|\cdot|z|^2(1+|x|)^{-d-\alpha-3}\stackrel{(\ref{ER60})}{\preceq}|x-x'|\cdot|z|^2p_\alpha(1,x).\label{ER8}
\end{align}
If $|z|>1$ and $|x-x'|\leq 1$, then we have
\begin{align}
|\delta_{p_\alpha}(1,x;z)-\delta_{p_\alpha}(1,x';z)|
&\preceq |x-x'|\int^1_0|\nabla p_\alpha(1,x\pm z+\theta(x'-x))|\dif\theta\no\\
&\quad+|x-x'|\int^1_0|\nabla p_\alpha(1,x+\theta(x'-x))|\dif\theta\no\\
&\stackrel{(\ref{ER5})}{\preceq} |x-x'|\left((1+|x\pm z|)^{-d-\alpha-1}+(1+|x|)^{-d-\alpha-1}\right)\no\\
&\stackrel{(\ref{ER60})}{\preceq}|x-x'|\left( p_\alpha(1,x\pm z)+p_\alpha(1,x)\right). \label{ER10}
\end{align}
Combining (\ref{ER9}), (\ref{ER8}) and (\ref{ER10}), we obtain
\begin{align*}
&|\delta_{p_\alpha}(1,x;z)-\delta_{p_\alpha}(1,x';z)|\\
&\quad\preceq (|x-x'|)\wedge 1)(|z|^2\wedge 1)
\left( p_\alpha(1,x\pm z)+p_\alpha(1,x)+p_\alpha(1,x'\pm z)+p_\alpha(1,x')\right),
\end{align*}
which implies (\ref{ER3}). The proof is complete.
\end{proof}

\subsection{Fractional derivative estimate of heat kernel of $\sL^\kappa_\alpha $}
Let $\kappa(z)$ be a measurable function on $\mR^d$ with
\begin{align}
\kappa(z)=\kappa(-z),\ \ 0<\kappa_0\leq \kappa(z)\leq \kappa_1.
\label{Ker}
\end{align}

Consider the following nonlocal symmetric operator
$$
\sL^\kappa_\alpha f(x):=\mbox{p.v.}\int_{\mR^d} (f(x+z)-f(x))\kappa(z)|z|^{-d-\alpha}\dif z=\frac{1}{2}\int_{\mR^d} \delta_f(x;z)\kappa(z)|z|^{-d-\alpha}\dif z,
$$
where $\delta_f(x;z)$ is defined in a similar way as in \eqref{e:2.10}
but with function $f$ not containing $t$ variable.
It is the infinitesimal generator of a symmetric L\'evy process that is
stable-like.
Let $p^\kappa_{\alpha} (t,x)$ be the heat kernel of operator $\sL^\kappa_\alpha $, i.e.,
$$
\p_t p^\kappa_{\alpha} (t,x)=\sL^\kappa_\alpha  p^\kappa_{\alpha} (t,x),\ \ \lim_{t\downarrow 0}p^\kappa_{\alpha} (t,x)=\delta_0(x).
$$
Under (\ref{Ker}), it is well-known
from the inverse Fourier transform that
\begin{align}
p^\kappa_\alpha\in C(\mR_+; C^\infty_b(\mR^d)).
\end{align}
Moreover, it follows from \cite[Theorem 1.1]{Ch-Ku} that
\begin{align}
p^\kappa_{\alpha} (t,x)\asymp \varrho^0_\alpha(t,x)=t(t^{1/\alpha}+|x|)^{-d-\alpha}.\label{ER66}
\end{align}
If we set
$$
\hat \kappa(z):=\kappa(z)-\tfrac{\kappa_0}{2},
$$
then by the construction of the L\'evy process,
one can write
\begin{align}
p^\kappa_{\alpha} (t,x)=\int_{\mR^d}p^{\kappa_0/2}_\alpha(t,x-y)p^{\hat \kappa}_\alpha(t,y)\dif y=\int_{\mR^d}p_\alpha\big(\tfrac{\kappa_0t}{2},x-y\big)p^{\hat\kappa}_\alpha(t,y)\dif y.\label{ER11}
\end{align}
The following lemma is an easy consequence of (\ref{ER66}), (\ref{ER11}) and Lemma \ref{Le4}.
\bl
Under (\ref{Ker}), there exists a constant $C=C(d,\alpha,\kappa_0,\kappa_1,
\kappa_2)>0$ such that
\begin{align}
|p^\kappa_{\alpha} (t,x)-p^\kappa_{\alpha} (t,x')|\leq C((t^{-1/\alpha}|x-x'|)\wedge 1)\left(\varrho^0_\alpha(t,x)+\varrho^0_\alpha(t,x')\right),\label{ER222}
\end{align}
\begin{align}
|\nabla p^\kappa_{\alpha}(t,x)|\leq Ct^{-1/\alpha}\varrho^0_\alpha(t,x),\label{ER70}
\end{align}
\begin{align}
|\delta_{p^\kappa_{\alpha} }(t,x;z)|\leq C
\left( (t^{-\frac{2}{\alpha}}|z|^2)\wedge 1 \right)
\left( \varrho^0_\alpha(t,x\pm z)+\varrho^0_\alpha(t,x)\right),\label{ER22}
\end{align}
\begin{align}
 |\delta_{p^\kappa_{\alpha} }(t,x;z)-\delta_{p^\kappa_{\alpha} }(t,x';z)|
& \leq C \left((t^{-1/\alpha}|x-x'|)\wedge 1\right) \left((t^{-\frac{2}{\alpha}}|z|^2)\wedge 1\right)\no\\
&\hskip 0.2truein \times\left( \varrho^0_\alpha(t,x\pm z)+\varrho^0_\alpha(t,x)+\varrho^0_\alpha(t,x'\pm z)+\varrho^0_\alpha(t,x')\right) .\label{ER33}
\end{align}
\el
\begin{proof}
By (\ref{ER11}) and (\ref{ER20}), we have
\begin{align*}
|p^\kappa_{\alpha} (t,x)-p^\kappa_{\alpha} (t,x')|
&\preceq ((t^{-1/\alpha}|x-x'|)\wedge 1)
\int_{\mR^d}\Big\{p_\alpha\big(\tfrac{\kappa_0t}{2},x-y\big)
+p_\alpha\big(\tfrac{\kappa_0t}{2},x'-y\big)\Big\}
p^{\hat \kappa}_\alpha(t,y)\dif y\\
&=((t^{-1/\alpha}|x-x'|)\wedge 1)\Big\{p^\kappa_{\alpha} (t,x)+p^\kappa_{\alpha} (t,x')\Big\}\\
&\stackrel{(\ref{ER66})}{\preceq}((t^{-1/\alpha}|x-x'|)\wedge 1)\Big\{\varrho^0_\alpha(t,x)+\varrho^0_\alpha(t,x')\Big\}.
\end{align*}
Similarly, we have (\ref{ER70}), (\ref{ER22}) and (\ref{ER33}) by (\ref{ER11}), (\ref{ER00}), (\ref{ER2}), (\ref{ER3}) and (\ref{ER66}).
\end{proof}
Now, we can prove the following fractional derivative estimate of $p^\kappa_\alpha(t,x)$.
\bt\label{Th24}
Under (\ref{Ker}), there exists a constant $C=C(d,\alpha,\kappa_0,\kappa_1,
\kappa_2)>0$ such that
\begin{align}
\int_{\mR^d}|\delta_{p^\kappa_{\alpha} }(t,x;z)|\cdot|z|^{-d-\alpha}\dif z&\leq C\varrho^0_0(t,x),\label{ER6}
\end{align}
\begin{align}
\int_{\mR^d}|\delta_{p^\kappa_{\alpha} }(t,x;z)-\delta_{p^\kappa_{\alpha} }(t,x';z)|\cdot |z|^{-d-\alpha}\dif z&\leq C((t^{-1/\alpha}|x-x'|)\wedge 1)\Big\{\varrho^0_0(t,x)+\varrho^0_0(t,x')\Big\}.\label{ER7}
\end{align}
\et
\begin{proof}
By (\ref{ER22}), we have
\begin{align*}
\int_{\mR^d}|\delta_{p^\kappa_{\alpha} }(t,x;z)|\cdot|z|^{-d-\alpha}\dif z
&\preceq\int_{\mR^d}((t^{-\frac{2}{\alpha}}|z|^2)\wedge 1)\varrho^0_\alpha(t,x\pm z)|z|^{-d-\alpha}\dif z\\
& \hskip 0.2truein +\varrho^0_\alpha(t,x)\int_{\mR^d}((t^{-\frac{2}{\alpha}}|z|^2)\wedge 1)|z|^{-d-\alpha}\dif z=:I_1+I_2.
\end{align*}
For $I_1$, we have
\begin{align*}
I_1&\leq t^{-\frac{2}{\alpha}}\int_{|z|\leq t^{1/\alpha}}\varrho^0_\alpha(t,x\pm z)|z|^{2-d-\alpha}\dif z
+\int_{|z|>t^{1/\alpha}}\varrho^0_\alpha(t,x\pm z)|z|^{-d-\alpha}\dif z=:I_{11}+I_{12}.
\end{align*}
For $I_{11}$, by (\ref{ER5}), we have
\begin{align*}
I_{11}&\preceq t^{1-\frac{2}{\alpha}}\int_{|z|\leq t^{1/\alpha}}(t^{1/\alpha}+|x\pm z|)^{-d-\alpha}|z|^{2-d-\alpha}\dif z\\
&\preceq t^{1-\frac{2}{\alpha}}(t^{1/\alpha}+|x|)^{-d-\alpha}\int_{|z|\leq t^{1/\alpha}}|z|^{2-d-\alpha}\dif z\preceq \varrho^0_0(t,x).
\end{align*}
For $I_{12}$, if $|x|\leq 2t^{1/\alpha}$, then
\begin{align*}
I_{12}&\preceq t\int_{|z|>t^{1/\alpha}}(t^{1/\alpha}+|x\pm z|)^{-d-\alpha}|z|^{-d-\alpha}\dif z\\
&\preceq t^{-d/\alpha}\int_{|z|>t^{1/\alpha}}|z|^{-d-\alpha}\dif z\preceq t^{-\frac{d+\alpha}{\alpha}}\leq\varrho^0_0(t,x);
\end{align*}
if $|x|>2t^{1/\alpha}$, then
\begin{align*}
I_{12}&\preceq \left(\int_{\frac{|x|}{2}\geq|z|>t^{1/\alpha}}+\int_{|z|>\frac{|x|}{2}}\right)\varrho^0_\alpha(t,x\pm z)\cdot|z|^{-d-\alpha}\dif z\\
&\preceq t\int_{\frac{|x|}{2}\geq|z|>t^{1/\alpha}}
(t^{1/\alpha}+|x\pm z|)^{-d-\alpha}|z|^{-d-\alpha}\dif z+|x|^{-d-\alpha}\int_{|z|>\frac{|x|}{2}}\varrho^0_\alpha(t,x\pm z)\dif z\\
&\preceq t(t^{1/\alpha}+|x|)^{-d-\alpha}\int_{|z|>t^{1/\alpha}}
|z|^{-d-\alpha}\dif z+|x|^{-d-\alpha}\int_{\mR^d}\varrho^0_\alpha(t,x\pm z)\dif z\\
&\preceq(t^{1/\alpha}+|x|)^{-d-\alpha}+|x|^{-d-\alpha}\preceq\varrho^0_0(t,x).
\end{align*}
For $I_2$, we have
\begin{align*}
I_2= t^{-1}\varrho^0_\alpha(t,x)\int_{\mR^d}(|z|^2\wedge 1)|z|^{-d-\alpha}\dif z\preceq\varrho^0_0(t,x).
\end{align*}
Combining the above calculations, we obtain (\ref{ER6}).

By (\ref{ER33}), as above, we have
\begin{align*}
&\int_{\mR^d}|\delta_{p^\kappa_{\alpha} }(t,x;z)-\delta_{p^\kappa_{\alpha} }(t,x';z)|\cdot |z|^{-d-\alpha}\dif z\preceq ((t^{-1/\alpha}|x-x'|)\wedge 1)\\
&\quad\times\Bigg\{\int_{\mR^d}((t^{-\frac{2}{\alpha}}|z|^2)\wedge 1)\{\varrho^0_\alpha(t,x\pm z)+\varrho^0_\alpha(t,x'\pm z)\}|z|^{-d-\alpha}\dif z\\
&\qquad+\{\varrho^0_\alpha(t,x)+\varrho^0_\alpha(t,x')\}\int_{\mR^d}((t^{-\frac{2}{\alpha}}|z|^2)\wedge 1)|z|^{-d-\alpha}\dif z\Bigg\}\\
&\quad\preceq ((t^{-1/\alpha}|x-x'|)\wedge 1)\Big\{\varrho^0_0(t,x)+\varrho^0_0(t,x')\Big\}.
\end{align*}
The proof is complete.
\end{proof}

\subsection{Continuous dependence of heat kernels with respect to $\kappa$}
In this subsection, we prove the following continuous dependence of
the heat kernel with respect to the kernel function $\kappa$, which seems to be new.

\bt \label{T:2.5}
Let $\kappa$ and $\tilde\kappa$ be two functions on $\mR^d$ satisfying (\ref{Ker}). For any $\gamma\in(0,\alpha\wedge 1)$, there exists a constant $C=C(d,\alpha,\kappa_0,\kappa_1, \kappa_2,\gamma)>0$ such that
\begin{align}
|p^{\kappa}_\alpha(t,x)-p^{\tilde\kappa}_\alpha(t,x)|\leq C\|\kappa-\tilde\kappa\|_\infty(\varrho^0_{\alpha}+\varrho^\gamma_{\alpha-\gamma})(t,x),\label{ER30}
\end{align}
\begin{align}
|\nabla p^{\kappa}_\alpha(t,x)-\nabla p^{\tilde\kappa}_\alpha(t,x)|\leq C\|\kappa-\tilde\kappa\|_\infty t^{-1/\alpha}(\varrho^0_{\alpha}+\varrho^\gamma_{\alpha-\gamma})(t,x),\label{ER308}
\end{align}
and
\begin{align}
\int_{\mR^d}|\delta_{p^{\kappa}_\alpha}(t,x;z)-\delta_{p^{\tilde\kappa}_\alpha}(t,x;z)|\cdot |z|^{-d-\alpha}\dif z\leq C\|\kappa-\tilde\kappa\|_\infty(\varrho^0_0+\varrho^\gamma_{-\gamma})(t,x).\label{ER77}
\end{align}
\et

\begin{proof}
(i) Note that the heat kernel
$p^\kappa_\alpha (t, x)$ is an even function in $x$.
We have
\begin{align*}
p^{\kappa}_\alpha(t,x)-p^{\tilde\kappa}_\alpha(t,x)
&=\int^t_0 \frac{\dif }{\dif s} \left( \int_{\mR^d}
p^{\kappa}_\alpha(s,y) p^{\tilde\kappa}_\alpha(t-s,x-y) \dif y\right) \dif s  \\
&=\int^t_0 \left(\int_{\mR^d} \left(\sL^{\kappa}_\alpha p^{\kappa}_\alpha(s,y) p^{\tilde\kappa}_\alpha(t-s,x-y)
-p^{\kappa}_\alpha(s,y) \sL^{\tilde\kappa}_\alpha p^{\tilde\kappa}_\alpha(t-s,x-y) \right)\dif y \right) \dif s\\
&=\int^t_0 \left( \int_{\mR^d}p^{\kappa}_\alpha(s,y)(\sL^{\kappa}_\alpha-\sL^{\tilde\kappa}_\alpha)p^{\tilde\kappa}_\alpha(t-s,x-y)\dif y \right) \dif s\\
&=\int^t_0 \left( \int_{\mR^d}(\sL^{\kappa}_\alpha-\sL^{\tilde\kappa}_\alpha)
p^{\tilde\kappa}_\alpha(t-s,x-y)\left(p^{\kappa}_\alpha(s,y)-
p^{\kappa}_\alpha(s,x)\right)\dif y \right) \dif s,
\end{align*}
where the third equality is due to the symmetry  of the
operator $\sL^{\tilde\kappa}_\alpha$, \eqref{ER66}, \eqref{ER6} and \eqref{EU7},
and the fourth equality is due to
$$
\int_{\mR^d}p^{\tilde\kappa}_\alpha(t-s,x-y)\dif y=1.
$$
Thus, by (\ref{ER222}) and (\ref{ER6}), we have
\begin{align*}
|p^{\kappa}_\alpha(t,x)-p^{\tilde\kappa}_\alpha(t,x)|&\leq\|\kappa-\tilde\kappa\|_\infty\int^t_0\!\!\!\int_{\mR^d}
\left(\int_{\mR^d}|\delta_{p^{\tilde\kappa}_\alpha}(t-s,x-y;z)|\cdot |z|^{-d-\alpha}\dif z\right)\\
&\hskip 1.2truein \times|p^{\kappa}_\alpha(s,y)-p^{\kappa}_\alpha(s,x)|\dif y\dif s\\
&\preceq\|\kappa-\tilde\kappa\|_\infty\int^t_0\!\!\!\int_{\mR^d}\varrho^0_0(t-s,x-y)\\
&\hskip 1.2truein  \times((s^{-1/\alpha}|x-y|)\wedge 1)(\varrho^0_\alpha(s,y)+\varrho^0_\alpha(s,x))\dif y\dif s\\
&\leq\|\kappa-\tilde\kappa\|_\infty\int^t_0\!\!\!\int_{\mR^d}\varrho^0_0(t-s,x-y)\\
&\hskip 1.2truein  \times((s^{-1/\alpha}|x-y|)^\gamma\wedge 1)(\varrho^0_\alpha(s,y)+\varrho^0_\alpha(s,x))\dif y\dif s\\
&\leq\|\kappa-\tilde\kappa\|_\infty\int^t_0\!\!\!\int_{\mR^d}\varrho^\gamma_0(t-s,x-y)(\varrho^0_{\alpha-\gamma}(s,y)+\varrho^0_{\alpha-\gamma}(s,x))\dif y\dif s\\
&\stackrel{(\ref{eq30})}{\preceq}\|\kappa-\tilde\kappa\|_\infty\Big\{\varrho^0_{\alpha}(t,x)+\varrho^\gamma_{\alpha-\gamma}(t,x)\Big\},
\end{align*}
which gives (\ref{ER30}).

(ii) By (\ref{ER11}), (\ref{ER70}) and (\ref{ER30}), we have
\begin{align*}
|\nabla p^{\kappa}_\alpha(t,x)-\nabla p^{\tilde\kappa}_\alpha(t,x)|&=\left|\int_{\mR^d}\nabla p_\alpha(\tfrac{\kappa_0t}{2},x-y)(p_{\hat \kappa}(t,y)-p_{\hat{ \tilde{\kappa}}}(t,y))\dif y\right|\\
&\preceq\|\kappa-\tilde\kappa\|_\infty t^{-1/\alpha}\int_{\mR^d}\varrho^0_\alpha(t,x-y)(\varrho^0_{\alpha}+\varrho^\gamma_{\alpha-\gamma})(t,y)\dif y\\
&\stackrel{(\ref{EU7})}{\preceq}\|\kappa-\tilde\kappa\|_\infty t^{-1/\alpha}(\varrho^0_{\alpha}+\varrho^\gamma_{\alpha-\gamma})(2t,x),
\end{align*}
which gives (\ref{ER308}).

(iii) By (\ref{ER11}), (\ref{ER22}) and (\ref{ER30}), we have
\begin{align*}
|\delta_{p^{\kappa}_\alpha}(t,x;z)-\delta_{p^{\tilde\kappa}_\alpha}(t,x;z)|&=\left|\int_{\mR^d}\delta_{p_\alpha}(\tfrac{\kappa_0t}{2},x-y;z)(p_{\hat \kappa}(t,y)-p_{\hat{ \tilde{\kappa}}}(t,y))\dif y\right|\\
&\preceq\|\kappa-\tilde\kappa\|_\infty((t^{-\frac{2}{\alpha}}|z|^2)\wedge 1)\\
&\times\int_{\mR^d}\Big\{\varrho^0_\alpha(t,x-y\pm z)+\varrho^0_\alpha(t,x-y)\Big\}(\varrho^0_{\alpha}+\varrho^\gamma_{\alpha-\gamma})(t,y)\dif y\\
&\stackrel{(\ref{EU7})}{\preceq}\|\kappa-\tilde\kappa\|_\infty((t^{-\frac{2}{\alpha}}|z|^2)\wedge 1)\\
&\times\Big\{(\varrho^0_{\alpha}+\varrho^\gamma_{\alpha-\gamma})(2t,x\pm z)+(\varrho^0_{\alpha}+\varrho^\gamma_{\alpha-\gamma})(2t,x)\Big\}.
\end{align*}
Using the same argument as in estimating (\ref{ER7}), we obtain (\ref{ER77}).
\end{proof}

\section{Levi's construction of heat kernels}

In this section we consider the spatial dependent operator
$\sL^{\kappa}_{\alpha}$ defined by \eqref{Non}, with the kernel
function $\kappa (x, z)$ satisfying conditions \eqref{Con1}-\eqref{Con2}.
 In order to reflect the dependence on $x$, we also write
$$
\sL^{\kappa(x)}_{\alpha} f(x)=\sL^{\kappa}_{\alpha} f(x)=\frac{1}{2}\int_{\mR^d} \delta_f(x;z)\kappa(x,z)|z|^{-d-\alpha}\dif z.
$$
For fixed $y\in\mR^d$, let $\sL^{\kappa(y)}_{\alpha}$ be the freezing operator
$$
\sL^{\kappa(y)}_{\alpha} f(x)=\frac{1}{2}\int_{\mR^d} \delta_f(x;z)\kappa(y,z)|z|^{-d-\alpha}\dif z.
$$
Let $p_y(t,x):=p^{\kappa(y)}_\alpha(t,x)$ be the heat kernel of operator $\sL^{\kappa(y)}_{\alpha}$, i.e.,
\begin{align}
\p_t p_y(t,x)=\sL^{\kappa(y)}_{\alpha} p_y(t,x),\ \ \lim_{t\downarrow 0}p_y(t,x)=\delta_0(x),\label{ES2}
\end{align}
where, with a little abuse of notation, $\delta_0(x)$ denotes the usual Dirac function.

Now, we want to seek the heat kernel $p^\kappa_{\alpha} (t,x,y)$ of $\sL^{\kappa}_{\alpha}$ with the following form:
\begin{align}
p^\kappa_{\alpha} (t,x,y)=p_y(t,x-y)+\int^t_0\!\!\!\int_{\mR^d}p_z(t-s,x-z)q(s,z,y)\dif z\dif s.\label{ER65}
\end{align}
The classical Levi's method
suggests that $q(t,x,y)$ solves the following integral equation:
\begin{align}
q(t,x,y)=q_0(t,x,y)+\int^t_0\!\!\!\int_{\mR^d}q_0(t-s,x,z)q(s,z,y)\dif z\dif s,\label{EU2}
\end{align}
where
$$
q_0(t,x,y):=(\sL^{\kappa(x)}_{\alpha}-\sL^{\kappa(y)}_{\alpha})p_y(t,x-y)=\int_{\mR^d}\delta_{p_y}(t,x-y;z)(\kappa(x,z)-\kappa(y,z))|z|^{-d-\alpha}\dif z.
$$
In fact, we formally have
\begin{align}
\p_tp^\kappa_{\alpha} (t,x,y)&=\sL^{\kappa(y)}_{\alpha}p_y(t,x-y)+q(t,x,y)+\int^t_0\!\!\!\int_{\mR^d}\p_tp_z(t-s,x-z)q(s,z,y)\dif z\dif s\no\\
&=\sL^{\kappa(x)}_{\alpha}p_y(t,x-y)+\int^t_0\!\!\!\int_{\mR^d}\sL^{\kappa(x)}_{\alpha}p_z(t-s,x-z)q(s,z,y)\dif z\dif s\no\\
&=\sL^{\kappa(x)}_{\alpha}p^\kappa_{\alpha} (t,x,y).\label{ER91}
\end{align}
Thus, the main aims of this section are to solve equation (\ref{EU2}), and to make the calculations in (\ref{ER91}) rigorous.

\subsection{Solving equation (\ref{EU2})}
In this subsection, we use Picard's iteration to solve (\ref{EU2}).
\bt\label{T3.4}
For $n\in\mN$, define $q_n(t,x,y)$ recursively by
\begin{align}
q_n(t,x,y):=\int^t_0\!\!\!\int_{\mR^d}q_0(t-s,x,z)q_{n-1}(s,z,y)\dif z\dif s.\label{EU22}
\end{align}
Under (\ref{Con1}) and (\ref{Con2}),  the series $q(t,x,y):=\sum_{n=0}^{\infty}q_n(t,x,y)$ is absolutely convergent and solves
the integral equation (\ref{EU2}). Moreover,
$q(t,x,y)$ has the following estimates:
there is a constant $C_1=C_1 (d, \alpha, \beta, \kappa_0, \kappa_1, \kappa_2)>0$ so that
\begin{align}
|q(t,x,y)| \leq C_1
(\varrho^\beta_0+\varrho^0_\beta)(t,x-y), \label{eq3}
\end{align}
and for any $\gamma\in(0,\beta)$,
there is a constant $C_2=C_2 (d, \alpha, \beta, \gamma, \kappa_0, \kappa_1, \kappa_2)>0$ so that
\begin{align}
|q(t,x,y)-q(t,x',y)| \leq C_2
\left( |x-x'|^{\beta-\gamma}\wedge 1\right) \left( (\varrho^0_\gamma+\varrho^\beta_{\gamma-\beta})(t,x-y)
+(\varrho^0_\gamma+\varrho^\beta_{\gamma-\beta})(t,x'-y)\right) .
\label{eq4}
\end{align}
\et
\begin{proof}
Without loss of generality, we assume $\beta\in(0,\frac{\alpha}{4}]$.
We divide the proof into three steps.

(Step 1). First of all, by (\ref{Con1}), (\ref{Con2}) and (\ref{ER6}), we have
\begin{align}
|q_0(t,x,y)|&\preceq(|x-y|^\beta\wedge 1)\int_{\mR^d}|\delta_{p_y}(t,x-y;z)|\cdot|z|^{-d-\alpha}\dif z\no\\
&\preceq(|x-y|^\beta\wedge 1)\varrho^0_0(t,x-y)=\varrho^\beta_0(t,x-y).\label{EU5}
\end{align}
For $n=1$, by definition (\ref{EU22}) and (\ref{eq30}), there exits a constant $C_{d,\alpha}>0$ such that
\begin{align}
|q_1(t,x,y)|\leq  C_{d,\alpha}\cB(\beta,\beta)\Big\{\varrho^0_{2\beta}+\varrho^\beta_\beta\Big\}(t,x-y).
\end{align}
Suppose now that
$$
|q_n(t,x,y)|\leq\gamma_n\Big\{\varrho^0_{(n+1)\beta}+\varrho^\beta_{n\beta}\Big\}(t,x-y),
$$
where $\gamma_n>0$ will be determined below. By (\ref{eq30}), we have
\begin{align*}
|q_{n+1}(t,x,y)|&\leq C_{d,\alpha}\gamma_{n}\cB(\beta,(n+1)\beta)\Big\{\varrho^0_{(n+2)\beta}+\varrho^\beta_{(n+1)\beta}\Big\}(t,x-y)\\
&=:\gamma_{n+1}\Big\{\varrho^0_{(n+2)\beta}+\varrho^\beta_{(n+1)\beta}\Big\}(t,x-y),
\end{align*}
where
$$
\gamma_{n+1}=C_{d,\alpha}\gamma_{n} \cB(\beta,(n+1)\beta).
$$
Hence, by $\cB(\gamma,\beta)=\frac{\Gamma(\gamma)\Gamma(\beta)}{\Gamma(\gamma+\beta)}$, where $\Gamma$ is the usual Gamma function,
we obtain
\begin{align*}
\gamma_{n}= C_{d,\alpha}^{n+1}\cB(\beta,\beta)\cB(\beta,2\beta)\cdots\cB(\beta,n\beta)=
\frac{(C_{d,\alpha}\Gamma(\beta))^{n+1}}{\Gamma((n+1)\beta)}.
\end{align*}
Thus,
\begin{align}
|q_n(t,x,y)|\leq\frac{(C_{d,\alpha}\Gamma(\beta))^{n+1}}{\Gamma((n+1)\beta)} \Big\{\varrho^0_{(n+1)\beta}+\varrho^\beta_{n\beta}\Big\}(t,x-y),\label{EYU1}
\end{align}
which in turn implies that
\begin{align*}
\sum_{n=0}^{\infty}|q_n(t,x,y)|&\leq
\sum_{n=0}^\infty\frac{(C_{d,\alpha}\Gamma(\beta))^{n+1}}{\Gamma((n+1)\beta)}\left\{\varrho^0_{(n+1)\beta}+\varrho^\beta_{n\beta}\right\}(t,x-y)\\
&\leq\sum_{n=0}^\infty\frac{(C_{d,\alpha}\Gamma(\beta))^{n+1}}{\Gamma((n+1)\beta)}\Big\{\varrho^0_{\beta}+\varrho^\beta_0\Big\}(t,x-y)\\
&\preceq \Big\{\varrho^0_{\beta}+\varrho^\beta_0\Big\}(t,x-y).
\end{align*}
Thus, (\ref{eq3}) is proven. Moreover, by (\ref{EU22}), we have
$$
\sum_{n=0}^{m+1}q_n(t,x,y)=q_0(t,x,y)+\int^t_0\!\!\!\int_{\mR^d}q_0(t-s,x,z)\sum_{n=0}^mq_n(s,z,y)\dif z\dif s,
$$
which yields (\ref{EU2}) by taking limits $m\to\infty$ for both sides.

(Step 2). In this step, we prove the following estimate:
\begin{align}
|q_0(t,x,y)-q_0(t,x',y)|
\preceq(|x-x'|^{\beta-\gamma}\wedge 1)
\Big\{(\varrho^0_\gamma+\varrho^\beta_{\gamma-\beta})(t,x-y)
+(\varrho^0_\gamma+\varrho^\beta_{\gamma-\beta})(t,x'-y)\Big\}.\label{ES5}
\end{align}
In the case of $|x-x'|>1$, we have
$$
|q_0(t,x,y)|\preceq \varrho^\beta_0(t,x-y)\leq\varrho^\beta_{\gamma-\beta}(t,x-y)
$$
and
$$
|q_0(t,x',y)|\preceq\varrho^\beta_0(t,x'-y)\leq\varrho^\beta_{\gamma-\beta}(t,x'-y).
$$
In the case of $1\geq |x-x'|>t^{1/\alpha}$, by (\ref{EU5}), we have
$$
|q_0(t,x,y)|\preceq\varrho^\beta_0(t,x-y)=t^{\frac{\beta-\gamma}{\alpha}}\varrho^\beta_{\gamma-\beta}(t,x-y)
\leq|x-x'|^{\beta-\gamma}\varrho^\beta_{\gamma-\beta}(t,x-y),
$$
and also
$$
|q_0(t,x',y)|\preceq|x-x'|^{\beta-\gamma}\varrho^\beta_{\gamma-\beta}(t,x'-y).
$$
Suppose now that
\begin{align}
|x-x'|\leq t^{1/\alpha}.\label{EE1}
\end{align}
By definition and Theorem \ref{Th24}, we have
\begin{align*}
|q_0(t,x,y)-q_0(t,x',y)|&=\Bigg|\int_{\mR^d}\delta_{p_y}(t,x-y,z)(\kappa(x,z)-\kappa(y,z))|z|^{-d-\alpha}\dif z\\
&\quad-\int_{\mR^d}\delta_{p_y}(t,x'-y,z)(\kappa(x',z)-\kappa(y,z))|z|^{-d-\alpha}\dif z\Bigg|\\
&\preceq(|x-y|^\beta\wedge 1)\int_{\mR^d}|\delta_{p_y}(t,x-y,z)-\delta_{p_y}(t,x'-y,z)|\cdot |z|^{-d-\alpha}\dif z\\
&\quad+(|x-x'|^\beta\wedge 1)\int_{\mR^d}|\delta_{p_y}(t,x'-y,z)|\cdot|z|^{-d-\alpha}\dif z\\
&\preceq (|x-y|^\beta\wedge 1)t^{-1/\alpha}|x-x'|\Big\{\varrho^0_0(t,x-y)+\varrho^0_0(t,x'-y)\Big\}\\
&\quad+(|x-x'|^\beta\wedge 1)\varrho^0_0(t,x'-y)\\
&\stackrel{(\ref{ER5})}{\preceq} t^{-1/\alpha}|x-x'|\varrho^\beta_0(t,x-y)+(|x-x'|^\beta\wedge 1)\varrho^0_0(t,x'-y)\\
&\stackrel{(\ref{EE1})}{\preceq} |x-x'|^{\beta-\gamma}\varrho^\beta_{\gamma-\beta}(t,x-y)+|x-x'|^{\beta-\gamma}\varrho^0_\gamma(t,x'-y).
\end{align*}
Combining the above calculations, we obtain (\ref{ES5}).

(Step 3). By definition (\ref{EU22}) and (\ref{EYU1}), (\ref{ES5}), we have for $n\in\mN$,
\begin{align*}
&|q_n(t,x,y)-q_n(t,x',y)|\preceq\int_0^t\!\!\!\int_{\mR^d}|q_0(t-s,x,z)-q_0(t-s,x',z)|q_{n-1}(s,z,y)\dif z\dif s\\
&\qquad\preceq\frac{(C_d\Gamma(\beta))^{n}}{\Gamma(n\beta)}\Big(|x-x'|^{\beta-\gamma}\wedge 1\Big)\int_0^t\!\!\!\int_{\mR^d}(\varrho_{n\beta}^{0}+\varrho_{(n-1)\beta}^{\beta})(s,z-y)\\
&\qquad\qquad\times\Big\{(\varrho_{\gamma}^0+\varrho_{\gamma-\beta}^{\beta})(t,x-z)+(\varrho_{\gamma}^0+\varrho_{\gamma-\beta}^{\beta})(t,x'-z)\Big\}\dif z\dif s\\
&\qquad\stackrel{(\ref{eq30})}{\preceq}\frac{(C_d\Gamma(\beta))^{n}}{\Gamma(n\beta)}\Big(|x-x'|^{\beta-\gamma}\wedge 1\Big)
\Big\{(\varrho^0_\gamma+\varrho^\beta_{\gamma-\beta})(t,x-y)+(\varrho^0_\gamma+\varrho^\beta_{\gamma-\beta})(t,x'-y)\Big\},
\end{align*}
which yields (\ref{eq4}) by summing up in $n$.
\end{proof}

\subsection{Some estimates about $p_y(t,x-y)$}
In this subsection, we prepare some important estimates for later use.
\bl
Under (\ref{Con1}) and (\ref{Con2}), there exists a constant
$C=C(d, \alpha, \beta, \kappa_0, \kappa_1, \kappa_2)>0$
 such that for all $\eps\geq 0$, $x,y\in\mR^d$ and $t>0$,
\begin{align}
\left|\int_{\mR^d}\left(\int_{|w|>\eps}\delta_{p_y}(t,x-y;w)\kappa(x,w)|w|^{-d-\alpha}\dif w\right)\dif y\right|\leq C t^{\frac{\beta}{\alpha}-1},\label{E303}
\end{align}
and
\begin{align}
\left|\int_{\mR^d}\nabla p_y(t,\cdot)(x-y)\dif y\right|\leq Ct^{\frac{\beta-1}{\alpha}}.\label{E300}
\end{align}
\el
\begin{proof}
Since
\begin{align}
\int_{\mR^d}p_x(t,\xi-y)\dif y=1, \  \ \forall\xi\in\mR^d,\label{EYU2}
\end{align}
by definition of $\delta_{p_x}(t,x-y;w)$, we have
$$
\int_{\mR^d}\delta_{p_x}(t, x-y;w)\dif y=0,\ \ \forall w\in\mR^d.
$$
Thus, by Fubini's theorem and (\ref{ER77}), we have for any $\gamma\in(0,\alpha\wedge 1)$,
\begin{align*}
&\left|\int_{\mR^d}\left(\int_{|w|>\eps}\delta_{p_y}(t,x-y;w)\kappa(x,w)|w|^{-d-\alpha}\dif w\right)\dif y\right|\\
&\quad=\left|\int_{\mR^d}\left(\int_{|w|>\eps}(\delta_{p_y}(t,x-y;w)-\delta_{p_x}(t,x-y;w))\kappa(x,w)|w|^{-d-\alpha}\dif w\right)\dif y\right|\\
&\quad\leq\kappa_1\int_{\mR^d}\left(\int_{|w|>\eps}|\delta_{p_y}(t,x-y;w)-\delta_{p_x}(t,x-y;w)|\cdot|w|^{-d-\alpha}\dif w\right)\dif y\\
&\quad\preceq\int_{\mR^d}\|\kappa(y,\cdot)-\kappa(x,\cdot)\|_\infty\Big\{\varrho^0_0(t,x-y)+\varrho^\gamma_{-\gamma}(t,x-y)\Big\}\dif y\\
&\quad\preceq\int_{\mR^d}\Big\{\varrho^\beta_0(t,x-y)+\varrho^{\beta+\gamma}_{-\gamma}(t,x-y)\Big\}\dif y\stackrel{(\ref{ES4})}{\preceq} t^{\frac{\beta}{\alpha}-1},
\end{align*}
which gives (\ref{E303}).

As for (\ref{E300}), it is similar by (\ref{EYU2}) and (\ref{ER308}) that
\begin{align*}
&\left|\int_{\mR^d}\nabla p_y(t,\cdot)(x-y)\dif y\right|=\left|\int_{\mR^d}(\nabla p_y(t,\cdot)-\nabla p_x(t,\cdot))(x-y)\dif y\right|\\
&\quad\preceq t^{-1/\alpha}\int_{\mR^d}\|\kappa(y,\cdot)-\kappa(x,\cdot)\|_\infty\Big\{\varrho^0_\alpha(t,x-y)+\varrho^\gamma_{\alpha-\gamma}(t,x-y)\Big\}\dif y\\
&\quad\preceq t^{-1/\alpha}\int_{\mR^d}\Big\{\varrho^\beta_\alpha(t,x-y)+\varrho^{\beta+\gamma}_{\alpha-\gamma}(t,x-y)\Big\}\dif y\stackrel{(\ref{ES4})}{\preceq} t^{\frac{\beta-1}{\alpha}}.
\end{align*}
The proof is complete.
\end{proof}
\bl
Under (\ref{Con1}) and (\ref{Con2}),
there is a constant $C=C(d, \alpha, \beta, \kappa_0, \kappa_1, \kappa_2)>0$
so that
\begin{align}
\left|\int_{\mR^d}\sL^{\kappa(x)}_{\alpha}p_y(t,\cdot)(x-y)\dif y\right|
\leq C t^{\frac{\beta}{\alpha}-1},\label{E30}
\end{align}
\begin{align}
\left|\int_{\mR^d}\p_tp_y(t,x-y)\dif y\right|
\leq C t^{\frac{\beta}{\alpha}-1},\label{E31}
\end{align}
\begin{align}
\lim_{t\downarrow 0}\sup_{x\in\mR^d}\left|\int_{\mR^d}p_y(t,x-y)\dif y-1\right|=0.\label{E32}
\end{align}
\el
\begin{proof}
Estimate  (\ref{E30}) follows by (\ref{E303}).
For (\ref{E31}), by (\ref{ES2}) we have
\begin{align*}
&\left|\int_{\mR^d}\p_tp_y(t,x-y)\dif y\right|=\left|\int_{\mR^d}\sL^{\kappa(y)}_{\alpha}p_y(t,\cdot)(x-y)\dif y\right|\\
&\quad\leq\left|\int_{\mR^d}(\sL^{\kappa(x)}_{\alpha}-\sL^{\kappa(y)}_{\alpha})p_y(t,\cdot)(x-y)\dif y\right|+\left|\int_{\mR^d}\sL^{\kappa(x)}_{\alpha}p_y(t,\cdot)(x-y)\dif y\right|\\
&\quad\stackrel{(\ref{EU5})(\ref{E30})}{\preceq}\int_{\mR^d}\varrho^\beta_0(t,x-y)\dif y+t^{\frac{\beta}{\alpha}-1}\preceq t^{\frac{\beta}{\alpha}-1}.
\end{align*}
For (\ref{E32}), by (\ref{EYU2}), we have  for any $\gamma\in(0,\alpha\wedge 1)$,
\begin{align*}
&\sup_{x\in\mR^d}\left|\int_{\mR^d}p_y(t,x-y)\dif y-1\right|\leq\sup_{x\in\mR^d}\int_{\mR^d}|p_y(t,x-y)-p_x(t,x-y)|\dif y\\
&\quad\stackrel{(\ref{ER30})}{\preceq}\sup_{x\in\mR^d}\int_{\mR^d}\|\kappa(y,\cdot)-\kappa(x,\cdot)\|_\infty(\varrho_\alpha(t,x-y)+\varrho^\gamma_{\alpha-\gamma}(t,x-y))\dif y\\
&\quad\preceq\sup_{x\in\mR^d}\int_{\mR^d}(\varrho^\beta_\alpha(t,x-y)+\varrho^{\gamma+\beta}_{\alpha-\gamma}(t,x-y))\dif y\preceq t^{\frac{\beta}{\alpha}}\to 0,
\end{align*}
as $t\to 0$. The proof is complete.
\end{proof}

\subsection{Smoothness of $p^\kappa_{\alpha} (t,x,y)$}
In this subsection, we give a rigorous proof about (\ref{ER91}).
Below, for the simplicity of notation, we write
\begin{align}
\phi_{y}(t,x,s):=\int_{\mR^d}p_z(t-s,x-z)q(s,z,y)\dif z,\label{ET77}
\end{align}
and
\begin{align}
\varphi_y(t,x):=\int^t_0\phi_y(t,x,s)\dif s=\int^t_0\!\!\!\int_{\mR^d}p_z(t-s,x-z)q(s,z,y)\dif z\dif s.\label{ET78}
\end{align}
First of all, we have
\bl\label{Le00}
For all $\gamma\in(0,\alpha\wedge 1)$,
there is a constant $C=C(d, \alpha, \beta, \gamma, \kappa_0, \kappa_1, \kappa_2)>0$
so that
$$
|p^\kappa_{\alpha} (t,x,y)-p^\kappa_{\alpha} (t,x',y)|
\leq C |x-x'|^\gamma\Big\{\varrho^0_{\alpha-\gamma}(t,x-y)
+\varrho^0_{\alpha-\gamma}(t,x'-y)\Big\}.
$$
\el
\begin{proof}
First of all, by (\ref{ER222}), we have
\begin{align*}
|p_y(t,x-y)-p_y(t,x'-y)|&\preceq((t^{-1/\alpha}|x-x'|)\wedge 1)\Big\{\varrho^0_\alpha(t,x-y)+\varrho^0_\alpha(t,x'-y)\Big\}\\
&\preceq|x-x'|^\gamma\Big\{\varrho^0_{\alpha-\gamma}(t,x-y)+\varrho^0_{\alpha-\gamma}(t,x'-y)\Big\}.
\end{align*}
On the other hand, by (\ref{eq3}) we also have
\begin{align*}
|\varphi_y(t,x)-\varphi_y(t,x')|&\leq\int^t_0\!\!\!\int_{\mR^d}|p_z(t-s,x-z)-p_z(t-s,x'-z)|\cdot|q(s,z,y)|\dif z\dif s\\
&\preceq\int^t_0\!\!\!\int_{\mR^d}((t-s)^{-1/\alpha}|x-x'|\wedge 1)\Big\{\varrho^0_\alpha(t-s,x-z)+\varrho^0_\alpha(t-s,x'-z)\Big\}\\
&\qquad\times\Big\{\varrho^0_\beta(s,z-y)+\varrho^\beta_0(s,z-y)\Big\}\dif z\dif s\\
&\preceq|x-x'|^\gamma\int^t_0\!\!\!\int_{\mR^d}\Big\{\varrho^0_{\alpha-\gamma}(t-s,x-z)+\varrho^0_{\alpha-\gamma}(t-s,x'-z)\Big\}\\
&\qquad\times\Big\{\varrho^0_\beta(s,z-y)+\varrho^\beta_0(s,z-y)\Big\}\dif z\dif s\\
&\preceq|x-x'|^\gamma\Big\{(\varrho^0_{\alpha-\gamma+\beta}+\varrho^\beta_{\alpha-\gamma})(t,x-y)+(\varrho^0_{\alpha-\gamma+\beta}+\varrho^\beta_{\alpha-\gamma})(t,x'-y)\Big\}.
\end{align*}
Combining the above two estimations, we obtain the desired estimate.
\end{proof}

\bl\label{Le44}
For all $x\not=y\in\mR^d$, the mapping $t\mapsto\varphi_y(t,x)$ is absolutely continuous, and
\begin{align}
\p_t \varphi_y(t,x)=q(t,x,y)+\int^t_0\!\!\!\int_{\mR^d}\sL^{\kappa(z)}_\alpha p_z(t-s,\cdot)(x-z)q(s,z,y)\dif z\dif s.\label{ES8}
\end{align}
\el
\begin{proof}
We divide the proof into four steps.

(Step 1). In this step we prove that for any $s\in(0,t)$,
\begin{align}
\p_t \phi_{y}(t,x,s)=\int_{\mR^d}\p_tp_z(t-s,x-z)q(s,z,y)\dif z.\label{EW1}
\end{align}
Notice that
\begin{align*}
\frac{\phi_{y}(t+\eps,x,s)-\phi_{y}(t,x,s)}{\eps}&=\frac{1}{\eps}\int_{\mR^d}\big(p_z(t+\eps-s,x-z)-p_z(t-s,x-z)\big)q(s,z,y)\dif z\no\\
&=\int_{\mR^d}\!\Bigg(\!\int_0^1 \p_t p_z(t+\theta \eps-s,x,z)\dif \theta \Bigg)q(s,z,y)\dif z.
\end{align*}
By (\ref{ES2}) and (\ref{ER6}), we have  for $|\eps|<\frac{t-s}{2}$,
\begin{align*}
|\p_t p_z(t+\theta \eps-s,x-z)|&=|\sL^{\kappa(z)}_\alpha p_z(t+\theta \eps-s,\cdot)(x-z)|\\
&\preceq (|x-z|+t+\theta\eps-s)^{-d-\alpha}\\
&\stackrel{(\ref{ER5})}{\preceq} (|x-z|+(t-s))^{-d-\alpha}\\
&=\varrho^0_0(t-s,x-z),
\end{align*}
which together with (\ref{eq3}) yields
\begin{align*}
|\p_t p_z(t+\theta \eps,x;r,z)q(s,z,y)|
\preceq \varrho_0^0(t-s,x-z)(\varrho^0_\beta+\varrho^\beta_0)(s,z-y)=:g(z).
\end{align*}
By (\ref{EU7}), one sees that
$$
\int_{\mR^d}g(z)\dif z<+\infty.
$$
Hence, by the dominated convergence theorem, we have
\begin{align*}
\lim_{\eps\to 0}\frac{\phi_{y}(t+\eps,x,s)-\phi_{y}(t,x,s)}{\eps}
=\int_{\mR^d}\p_tp_z(t-s,x-z)q(s,z,y)\dif z,
\end{align*}
and (\ref{EW1}) is proven.

(Step 2). In this step we prove that  for all $x\neq y$ and $t>0$,
\begin{align}
\int^t_0\!\!\!\int^r_0|\p_r \phi_{y}(r,x,s)|\dif s\dif r<+\infty.\label{EU9}
\end{align}
By (\ref{EW1}), we have
\begin{align}
|\p_r \phi_{y}(r,x,s)|&\leq\int_{\mR^d}|\p_rp_z(r-s,x-z)|\cdot|q(s,z,y)-q(s,x,y)|\dif z\no\\
&\quad+|q(s,x,y)|\left|\int_{\mR^d}\p_rp_z(r-s,x-z)\dif z\right|\no\\
&=:Q^{(1)}_y(r,x,s)+Q^{(2)}_y(r,x,s).\label{EP4}
\end{align}
For $Q^{(1)}_{y}(r,x,s)$, by (\ref{eq4}) and (\ref{ER6}), we have
\begin{align}
\int^t_0\!\!\!\int^r_0Q^{(1)}_{y}(r,x,s)\dif s\dif r&\preceq
\int^t_0\!\!\!\int^r_0\!\!\!\int_{\mR^d}|\sL^{\kappa(z)}_\alpha p_z(r-s,x-z)|\cdot (|x-z|^{\beta-\gamma}\wedge 1)\no\\
&\quad\times\Big\{(\varrho^0_\gamma+\varrho^\beta_{\gamma-\beta})(s,x-y)+
(\varrho^0_\gamma+\varrho^\beta_{\gamma-\beta})(s,z-y)\Big\}\dif z\dif s\dif r\no\\
&\preceq\int^t_0\!\!\!\int^r_0\!\!\!\int_{\mR^d}\varrho^{\beta-\gamma}_0(r-s,x-z)
(\varrho^0_\gamma+\varrho^\beta_{\gamma-\beta})(s,x-y)\dif z\dif s\dif r\no\\
&\quad+\int^t_0\!\!\!\int^r_0\!\!\!\int_{\mR^d}\varrho^{\beta-\gamma}_0(r-s,x-z)
(\varrho^0_\gamma+\varrho^\beta_{\gamma-\beta})(s,z-y)\dif z\dif s\dif r\no\\
&\preceq \int^t_0\!\!\!\int^r_0(r-s)^{\frac{\beta-\gamma}{\alpha}-1}(\varrho^0_\gamma+\varrho^\beta_{\gamma-\beta})(s,x-y)\dif s\dif r\no\\
&\quad+\int^t_0(\varrho^0_\beta+\varrho^\beta_0+\varrho^{\beta-\gamma}_\gamma)(r,x-y)\dif r\no\\
&\preceq \frac{1}{|x-y|^{d+\alpha}}\int^t_0\!\!\!\int^{r}_0(r-s)^{\frac{\beta-\gamma}{\alpha}-1}(s^{\frac{\gamma}{\alpha}}+s^{\frac{\gamma-\beta}{\alpha}})\dif s\dif r\no\\
&\quad+\frac{1}{|x-y|^{d+\alpha}}\int^t_0(r^{\frac{\gamma}{\alpha}}+1+r^{\frac{\beta}{\alpha}})\dif r<+\infty.\label{EP5}
\end{align}
For $Q^{(2)}_{y}(r,x,s)$, by (\ref{E31}) and (\ref{eq3}) we have
\begin{align}
\int^t_0\!\!\!\int^{r}_0Q^{(2)}_{y}(r,x,s)\dif r\dif r&\preceq\int^t_0\!\!\!\int^{r}_0
(\varrho^0_\beta+\varrho^\beta_0)(s,x-y)(r-s)^{\frac{\beta}{\alpha}-1}\dif s\dif r<+\infty.\label{EP6}
\end{align}
Combining (\ref{EP4})-(\ref{EP6}), we obtain (\ref{EU9}).

(Step 3). For fixed $s,x,y$, we have
\begin{align}
\lim_{t\downarrow s}\phi_{y}(t,x,s)=q(s,x,y).\label{EW11}
\end{align}
By (\ref{E32}), it suffices to prove that
$$
\lim_{t\downarrow r}\Bigg|\int_{\mR^d}p_z(t-s,x-z)(q(s,z,y)- q(s,x,y))\dif z \Bigg|=0.
$$
Notice that for any $\delta>0$,
\begin{align*}
&\Bigg|\int_{\mR^d}p_z(t-s,x-z)(q(s,z,y)- q(s,x,y))\dif z \Bigg|\\
&\quad \leq\int_{|x-z|\leq\delta}p_z(t-s,x-z)|q(s,z,y)-q(s,x,y)|\dif z\\
&\quad +\int_{|x-z|>\delta}p_z(t-s,x-z)|q(s,z,y)-q(s,x,y)|\dif z\\
&\quad =: J_1(\delta,t,s)+J_2(\delta,t,s).
\end{align*}
For any $\eps>0$, by (\ref{eq4}), there exists a $\delta=\delta(s,x,y)>0$ such that
for all $|x-z|\leq\delta$,
$$
|q(s,z,y)-q(s,x,y)|\leq \eps.
$$
Thus,
\begin{align*}
J_1(\delta,t,s)&\leq \varepsilon\int_{|x-z|\leq \delta}p_z(t-s,x-z)\dif z\\
&\leq\varepsilon\int_{\mR^d}p_z(t-s,x-z)\dif z\\
&\preceq\eps\int_{\mR^d}\varrho^0_\alpha(t-s,x-z)\dif z\stackrel{(\ref{ES4})}{\preceq}\eps.
\end{align*}
On the other hand, we have
\begin{align*}
J_2(\delta,t,s)&\stackrel{(\ref{ER66})}{\preceq}(t-s)\int_{|x-z|>\delta}\frac{|q(s,z,y)|+|q(s,x,y)|}{|x-z|^{d+\alpha}}\dif z\\
&\leq(t-s)\left(\delta^{-d-\alpha}\int_{\mR^d}|q(s,z,y)|\dif z
+|q(s,x,y)|\int_{|z|>\delta}|z|^{-d-\alpha}\dif z\right),
\end{align*}
which, by (\ref{eq3}) and (\ref{ES4}), converges to zero as $t\downarrow r$. Thus,  (\ref{EW11}) is proved.

(Step 4). Now, by the integration by parts formula and (\ref{EW11}), we have
$$
\int^t_s\p_{r} \phi_{y}(r,x,s)\dif r=\phi_{y}(t,x,s)-q(s,x,y).
$$
Integrating both sides with respect to $s$ from $0$ to $t$, and then by (\ref{EU9}) and Fubini's theorem, we obtain
\begin{align*}
\varphi_y(t,x)-\int^t_0q(s,x,y)\dif s&=\int^t_0\!\!\!\int^t_s\p_{r} \phi_{y}(r,x,s)\dif r\dif s
\stackrel{(\ref{EU9})}{=}\int^t_0\!\!\!\int^r_0\p_{r} \phi_{y}(r,x,s)\dif s\dif r\\
&\stackrel{(\ref{EW1}) (\ref{ES2})}{=}\int^t_0\!\!\!\int^r_0\!\!\!\int_{\mR^d}\sL^{\kappa(z)}_\alpha p_z(r-s,\cdot)(x-z)q(s,z,y)\dif z\dif s\dif r,
\end{align*}
which in turn implies (\ref{ES8}) by the Lebesgue differential theorem.
\end{proof}

\bl\label{Le6}
For all $t>0$ and $x\not=y$, we have
\begin{align}
\sL^{\kappa(x)}_{\alpha}\varphi_y(t,x)&=\int^t_0\!\!\!\int_{\mR^d}\sL^{\kappa(x)}_{\alpha}p_z(t-s,\cdot)(x-z)q(s,z,y)\dif z\dif s,\label{EP1}
\end{align}
and if $\beta>(1-\alpha)\vee 0$, then
\begin{align}
\nabla\varphi_y(t,x)&=\int^t_0\!\!\!\int_{\mR^d}\nabla p_z(t-s,\cdot)(x-z)q(s,z,y)\dif z\dif s,\label{EP11}
\end{align}
where the integrals are understood in the sense of iterated integrals. Moreover, for any $x\not=y$,
\begin{align}
\mbox{$t\mapsto \sL^{\kappa(x)}_{\alpha}\varphi_y(t,x)$ is continuous on $(0,1)$.}\label{EP22}
\end{align}
\el
\begin{proof}
We only prove (\ref{EP1}), and  (\ref{EP11}) is analogue by using (\ref{E300}).
First of all, for fixed $s\in(0,t)$, since
$$
x\mapsto p_y(t-s,x-y)\in C^\infty_b(\mR^d\times\mR^d)
$$
and
$$
z\mapsto q(s,z,y)\in C_b(\mR^d),
$$
by (\ref{ER6}) and Fubini's theorem, it is easy to see that
\begin{align}
\sL^{\kappa(x)}_{\alpha}\phi_{y}(t,x,s)=\int_{\mR^d}\sL^{\kappa(x)}_{\alpha} p_z(t-s,\cdot)(x-z)q(s,z,y)\dif z.\label{EY3}
\end{align}
By definition of $\phi_y$ and Fubini's theorem, we have for $\eps\in(0,1)$
\begin{align*}
I_\eps(t,x,s,y)&:=\left|\int_{|w|>\eps}\delta_{\phi_{y}}(t,x,s;w)\kappa(x,w)|w|^{-d-\alpha}\dif w\right|\\
&=\left|\int_{|w|>\eps}\left(\int_{\mR^d}\delta_{p_{z}}(t-s,x-z;w) q(s,z,y)\dif z\right) \kappa(x,w)|w|^{-d-\alpha}\dif w\right|\\
&=\left|\int_{\mR^d}\left(\int_{|w|>\eps}\delta_{p_{z}}(t-s,x-z;w) \kappa(x,w)|w|^{-d-\alpha}\dif w\right)q(s,z,y) \dif z\right|\\
&\leq\int_{\mR^d}\left(\int_{|w|>\eps}|\delta_{p_{z}}(t-s,x-z;w)|\cdot |w|^{-d-\alpha}\dif w\right)|q(s,z,y)-q(s,x,y)| \dif z\\
&\quad+\left|\int_{\mR^d}\left(\int_{|w|>\eps}\delta_{p_{z}}(t-s,x-z;w) \kappa(x,w)|w|^{-d-\alpha}\dif w\right) \dif z\right|\cdot|q(s,x,y)|.
\end{align*}
Using (\ref{ER6}), (\ref{E303}), (\ref{eq3}) and (\ref{eq4}), we further have
\begin{align*}
I_\eps(t,x,s,y)&\preceq \int_{\mR^d}\varrho^{\beta-\gamma}_0(t-s,x-z)(\varrho^0_\gamma+\varrho^\beta_{\gamma-\beta})(s,z-y)\dif z\\
&\quad+\left(\int_{\mR^d}\varrho^{\beta-\gamma}_0(t-s,x-z)\dif z\right)(\varrho^0_\gamma+\varrho^\beta_{\gamma-\beta})(s,x-y)\\
&\quad+(t-s)^{\frac{\beta}{\alpha}-1}(\varrho^\beta_0(s,x-y)+\varrho^0_\beta(s,x-y))\\
&\preceq \int_{\mR^d}\varrho^{\beta-\gamma}_0(t-s,x-z)(\varrho^0_\gamma+\varrho^\beta_{\gamma-\beta})(s,z-y)\dif z\\
&\quad+(t-s)^{\frac{\beta-\gamma}{\alpha}-1}(\varrho^0_\gamma+\varrho^\beta_{\gamma-\beta})(s,x-y)\\
&\quad+(t-s)^{\frac{\beta}{\alpha}-1}(\varrho^\beta_0(s,x-y)+\varrho^0_\beta(s,x-y)),
\end{align*}
which implies that for some $p>1$,
\begin{align}
\sup_{\eps\in(0,1)}\int^t_0|I_\eps(t,x,s,y)|^p\dif s<+\infty.\label{ETR1}
\end{align}
Now, by Fubini's theorem again, we obtain
\begin{align*}
\sL^{\kappa(x)}_{\alpha}\varphi_y(t,x)&=\lim_{\eps\downarrow 0}\int_{|w|>\eps}\!\int^t_0
\delta_{\phi_{y}}(t,x,s;w)\kappa(x,w)|w|^{-d-\alpha}\dif s\dif w\\
&=\lim_{\eps\downarrow 0}\int^t_0\!\!\!\int_{|w|>\eps}
\delta_{\phi_{y}}(t,x,s;w)\kappa(x,w)|w|^{-d-\alpha}\dif w\dif r\\
&=\int^t_0\lim_{\eps\downarrow 0}\int_{|w|>\eps}
\delta_{\phi_{y}}(t,x,s;w)\kappa(x,w)|w|^{-d-\alpha}\dif w\dif r\\
&=\int^t_0\sL^{\kappa(x)}_{\alpha}\phi_{y}(t,x,s)\dif s,
\end{align*}
which together with (\ref{EY3}) yields (\ref{EP1}).

As for (\ref{EP22}), it follows by (\ref{EP1}) and a direct calculation.
\end{proof}

\section{Proofs of Theorem \ref{Main} and Corollary \ref{C:1.4}}

\subsection{A nonlocal maximal principle}
In this subsection, we prove a nonlocal maximal principle (cf. \cite{Zh1}).
Notice that the current assumptions are weaker than \cite{Zh1}.
\bt\label{Th1}
Let $u(t,x)\in C_b([0,1]\times\mR^d)$ with
\begin{align}
\lim_{t\downarrow 0}\sup_{x\in\mR^d}|u(t,x)-u(0,x)|=0.\label{ER90}
\end{align}
Suppose that for each $x\in\mR^d$,
\begin{align}
\mbox{$t\mapsto \sL^\kappa_\alpha  u(t,x)$ is continuous on $(0,1]$},\label{ER93}
\end{align}
and for any $\eps\in(0,1)$ and some $\gamma\in((\alpha-1)\vee 0,1)$,
\begin{align}
\sup_{t\in(\eps,1)}|u(t,x)-u(t,x')|\leq C_\eps|x-x'|^\gamma.\label{ER92}
\end{align}
If $u(t,x)$  satisfies the following equation: for all $(t,x)\in(0,1)\times\mR^d$,
$$
\p_t u(t,x)=\sL^\kappa_\alpha  u(t,x),
$$
then for all $t\in(0,1)$,
$$
\sup_{x\in\mR^d}u(t,x)\leq \sup_{x\in\mR^d}u(0,x).
$$
\et
\begin{proof}
First of all, by (\ref{ER90}), it suffices to prove that for any $\eps\in(0,1)$,
\begin{align}
\sup_{x\in\mR^d}u(t,x)\leq \sup_{x\in\mR^d}u(\eps,x),\ \ \forall t\in(\eps,1).\label{ER911}
\end{align}
Below, we shall fix $\eps\in(0,1)$.
Let $\chi(x):\mR^d\to[0,1]$ be a smooth function with $\chi(x)=1$ for $|x|\leq 1$ and $\chi(x)=0$ for $|x|>2$. For $R>0$, define the following cutoff function
$$
\chi_R(x):=\chi(x/R).
$$
For $R,\delta>0$, consider
$$
u^\delta_R(t,x):=u(t,x)\chi_R(x)- (t-\eps)\delta.
$$
Then
\begin{align}
\p_t u^\delta_R(t,x)=\sL^\kappa_\alpha  u^\delta_R(t,x)+g^\delta_R(t,x),\label{EQ6}
\end{align}
where
$$
g^\delta_R(t,x):=\sL^\kappa_\alpha  u(t,x)\chi_R(x)-\sL^\kappa_\alpha  (u\chi_R)(t,x)-\delta.
$$
Our aim is to prove that for each $\delta>0$, there exists an $R_0\geq 1$ such that for all $t\in(\eps,1)$ and $R>R_0$,
\begin{align}
\sup_{x\in\mR^d}u^\delta_R(t,x) \leq \sup_{x\in\mR^d}u^\delta_R(\eps,x)\leq\sup_{x\in\mR^d}u(\eps,x). \label{EY1}
\end{align}
If this is proven, then taking $R\to\infty$ and $\delta\to 0$, we obtain (\ref{ER911}).

We first prove the following claim:

{\it Claim}:  For $\beta\in(0,\alpha\wedge 1)$, there exists a constant $C_\eps>0$ such that for all $R\geq 1$,
\begin{align}
\sup_{(t,x)\in[\eps,1]\times\mR^d}|\sL^\kappa_\alpha  u(t,x)\chi_R(x)-\sL^\kappa_\alpha  (u\chi_R)(t,x)|\leq \frac{C_\eps}{R^\beta}.\label{EQ4}
\end{align}
Moreover, for each $x\in\mR^d$,
\begin{align}
\mbox{$t\mapsto \sL^\kappa_\alpha u^\delta_R(t,x)$ and $g^\delta_R(t,x)$ are continuous on $(\eps,1)$}.\label{ER94}
\end{align}
{\it Proof of Claim:} Notice that by definitions,
\begin{align}
&\sL^\kappa_\alpha  (u\chi_R)(t,x)-\sL^\kappa_\alpha  u(t,x)\chi_R(x)-u(t,x)\sL^\kappa_\alpha \chi_R(x)\no\\
&\quad=\int_{\mR^d}(u(t,x+z)-u(t,x))(\chi_R(x+z)-\chi_R(x))\kappa(x,z)|z|^{-d-\alpha}\dif z.\label{EQ0}
\end{align}
Thus,
\begin{align*}
&|\sL^\kappa_\alpha  (u\chi_R)(t,x)-\sL^\kappa_\alpha  u(t,x)\chi_R(x)-u(t,x)\sL^\kappa_\alpha \chi_R(x)|\\
&\quad\leq\|\kappa\|_\infty\int_{|z|>1}|u(t,x+z)-u(t,x)|\cdot|\chi_R(x+z)-\chi_R(x)|\cdot|z|^{-d-\alpha}\dif z,\\
&\quad+\|\kappa\|_\infty\int_{|z|\leq 1}|u(t,x+z)-u(t,x)|\cdot|\chi_R(x+z)-\chi_R(x)|\cdot|z|^{-d-\alpha}\dif z=I_1+I_2.
\end{align*}
For $I_1$, we have
\begin{align}
I_1\leq 2\|\kappa\|_\infty\|u\|_\infty\int_{|z|>1}(2\|\chi_R\|_\infty)^{1-\beta}\|\chi'_R\|^{\beta}_\infty |z|^{\beta-d-\alpha}\dif z\preceq\|\kappa\|_\infty\|u\|_\infty(2\|\chi\|_\infty)^{1-\beta}\|\chi'\|^{\beta}_\infty/R^{\beta}.\label{EQ1}
\end{align}
For $I_2$, by (\ref{ER92}), we have
\begin{align}
I_2\leq \|\kappa\|_\infty C_\eps\int_{|z|\leq1}\|\chi'_R\|_\infty|z|^{1+\gamma-d-\alpha}\dif z\preceq \|\kappa\|_\infty C_\eps\|\chi'\|_\infty/R.\label{EQ2}
\end{align}
Moreover, it is also easy to see that
\begin{align}
\|\sL^\kappa_\alpha \chi_R\|_\infty\leq \frac{C}{R^{\beta}}.\label{EQ5}
\end{align}
Combining (\ref{EQ0})-(\ref{EQ5}), we obtain (\ref{EQ4}). As for (\ref{ER94}), it follows by (\ref{ER93}), (\ref{EQ0}) and the dominated convergence theorem.

\vspace{5mm}

We now use the contradiction argument to prove (\ref{EY1}). Fix
\begin{align}
R>(2C_\eps/\delta)^{1/\beta}.\label{ER95}
\end{align}
Suppose that (\ref{EY1}) does not hold, then there exists a $(t_0,x_0)\in(\eps,1)\times\mR^d$ such that
\begin{align}
\sup_{(t,x)\in(\eps,1)\times\mR^d}u^\delta_R(t,x)=u^\delta_R(t_0,x_0).\label{EYU3}
\end{align}
Thus, by (\ref{EQ6}), we have for any $h\in(0,t_0-\eps)$,
$$
0\leq \frac{u^\delta_R(t_0,x_0)-u^\delta_R(t_0-h,x_0)}{h}=\frac{1}{h}\int^{t_0}_{t_0-h}\sL^\kappa_\alpha  u^\delta_R(s,x_0)\dif s+\frac{1}{h}\int^{t_0}_{t_0-h}g^\delta_R(s,x_0)\dif s,
$$
which implies by  (\ref{ER94}) and letting $h\to 0$ that
\begin{align}
0\leq\sL^\kappa_\alpha  u^\delta_R(t_0,x_0)+g^\delta_R(t_0,x_0).\label{ER96}
\end{align}
On the other hand, by definition of $\sL^\kappa_\alpha $ and (\ref{EYU3}), we have
\begin{align}
\sL^\kappa_\alpha  u^\delta_R(t_0,x_0)=\int_{\mR^d}\delta_{u^\delta_R}(t_0,x_0;z)\kappa(x_0,z)|z|^{-d-\alpha}\dif z\leq 0,\label{ER97}
\end{align}
and by the claim and (\ref{ER95}),
\begin{align}
g^\delta_R(t_0,x_0)\leq\frac{C_\eps}{R^{\beta}}-\delta\leq-\frac{\delta}{2}.\label{ER98}
\end{align}
Combining (\ref{ER96})-(\ref{ER98}), we obtain a contradiction, and the proof is complete.
\end{proof}
\subsection{Fractional derivative and gradient estimates of $p^\kappa_\alpha$} We prove two lemmas about the fractional derivative and gradient estimates of $p^\kappa_\alpha$.
\bl\label{Le01}
We have
\begin{align}
|\sL^\kappa_\alpha p^\kappa_{\alpha} (t,\cdot,y)(x)|\preceq \varrho^0_0(t,x-y),\label{GH1}
\end{align}
and if $\alpha\in[1,2)$, then
\begin{align}
|\nabla p^\kappa_{\alpha} (t,x,y)|\preceq t^{\frac{\alpha-1}{\alpha}}\varrho^0_0(t,x-y).\label{GH2}
\end{align}
\el
\begin{proof}
(i) First of all, by (\ref{ER6}), it is easy to see that
$$
|\sL^{\kappa}_{\alpha}p_y(t,\cdot)(x-y)|\preceq \varrho^0_0(t,x-y).
$$
Recalling (\ref{ET78}), by (\ref{EP1}), we can write
\begin{align*}
\sL^{\kappa}_{\alpha}\varphi_y(t,x)&=\int^t_{\frac{t}{2}}\!\!\!\int_{\mR^d}\sL^{\kappa}_{\alpha}p_z(t-s,\cdot)(x-z)(q(s,z,y)-q(s,x,y))\dif z\dif s\\
&\quad+\int^t_{\frac{t}{2}}\left(\int_{\mR^d}\sL^{\kappa}_{\alpha}p_z(t-s,\cdot)(x-z)\dif z\right)q(s,x,y)\dif s\\
&\quad+\int^{\frac{t}{2}}_0\!\!\!\int_{\mR^d}\sL^{\kappa}_{\alpha}p_z(t-s,\cdot)(x-z)q(s,z,y)\dif z\dif s\\
&=:Q_1(t,x,y)+Q_2(t,x,y)+Q_3(t,x,y).
\end{align*}
For $Q_1(t,x,y)$, by (\ref{ER6}) and (\ref{eq4}), we have  for any $\gamma\in(0,\beta)$,
\begin{align*}
Q_1(t,x,y)&\preceq\int^t_{\frac{t}{2}}\!\!\!\int_{\mR^d}\varrho^{\beta-\gamma}_0(t-s,x-z)(\varrho^0_\gamma+\varrho^\beta_{\gamma-\beta})(s,x-y)\dif z\dif s\\
&+\int^t_{\frac{t}{2}}\!\!\!\int_{\mR^d}\varrho^{\beta-\gamma}_0(t-s,x-z)(\varrho^0_\gamma+\varrho^\beta_{\gamma-\beta})(s,z-y)\dif z\dif s\\
&\preceq\int^t_{\frac{t}{2}}(t-s)^{\frac{\beta-\gamma}{\alpha}-1}(\varrho^0_\gamma+\varrho^\beta_{\gamma-\beta})(s,x-y)\dif s\\
&+(\varrho^0_\beta+\varrho^{\beta-\gamma}_\gamma+\varrho^\beta_0)(t,x-y)\preceq\varrho^0_0(t,x-y).
\end{align*}
For $Q_2(t,x,y)$, by (\ref{E30}), we have
$$
Q_2(t,x,y)\preceq\int^t_{\frac{t}{2}}(t-s)^{\frac{\beta}{\alpha}-1}(\varrho^0_\beta+\varrho^\beta_0)(s,x-y)\dif s\preceq\varrho^0_0(t,x-y).
$$
For $Q_3(t,x,y)$, by  (\ref{ER6}), (\ref{eq3}) and (\ref{EU7}), we have
\begin{align*}
Q_3(t,x,y)&\preceq\int^{\frac{t}{2}}_0\!\!\!\int_{\mR^d}\varrho^0_0(t-s,x-z)(\varrho^0_\beta+\varrho^\beta_0)(s,z-y)\dif z\dif s\preceq\varrho^0_0(t,x-y).
\end{align*}
Combining the above calculations and by (\ref{ER65}), we obtain (\ref{GH1}).

(ii) By (\ref{ER70}), we have
$$
|\nabla p_y(t,\cdot)(x-y)|\preceq t^{\frac{\alpha-1}{\alpha}} \varrho^0_0(t,x-y).
$$
By (\ref{EP11}), we can write
\begin{align*}
\nabla\varphi_y(t,x)&=\int^t_{\frac{t}{2}}\!\!\!\int_{\mR^d}\nabla p_z(t-s,\cdot)(x-z)(q(s,z,y)-q(s,x,y))\dif z\dif s\\
&\quad+\int^t_{\frac{t}{2}}\left(\int_{\mR^d}\nabla p_z(t-s,\cdot)(x-z)\dif z\right)q(s,x,y)\dif s\\
&\quad+\int^{\frac{t}{2}}_0\!\!\!\int_{\mR^d}\nabla p_z(t-s,\cdot)(x-z)q(s,z,y)\dif z\dif s\\
&=:R_1(t,x,y)+R_2(t,x,y)+R_3(t,x,y).
\end{align*}
For $R_1(t,x,y)$, by (\ref{ER70}), (\ref{eq4}) and Lemma \ref{Le09},
in view of $\alpha\in[1,2)$, we have  for any $\gamma\in(0,\beta)$,
\begin{align*}
R_1(t,x,y)&\preceq\int^t_{\frac{t}{2}}\!\!\!\int_{\mR^d}\varrho^{\beta-\gamma}_{\alpha-1}(t-s,x-z)(\varrho^0_\gamma+\varrho^\beta_{\gamma-\beta})(s,x-y)\dif z\dif s\\
&+\int^t_{\frac{t}{2}}\!\!\!\int_{\mR^d}\varrho^{\beta-\gamma}_{\alpha-1}(t-s,x-z)(\varrho^0_\gamma+\varrho^\beta_{\gamma-\beta})(s,z-y)\dif z\dif s\\
&\preceq\int^t_{\frac{t}{2}}(t-s)^{\frac{\beta-\gamma-1}{\alpha}}(\varrho^0_\gamma+\varrho^\beta_{\gamma-\beta})(s,x-y)\dif s\\
&+(\varrho^0_{\beta+\alpha-1}+\varrho^{\beta-\gamma}_{\alpha+\gamma-1}+\varrho^\beta_{\alpha-1})(t,x-y)\preceq\varrho^0_{\alpha-1}(t,x-y).
\end{align*}
For $R_2(t,x,y)$, by (\ref{E300}), we have
$$
R_2(t,x,y)\preceq\int^t_{\frac{t}{2}}(t-s)^{\frac{\beta-1}{\alpha}}(\varrho^0_\beta+\varrho^\beta_0)(s,x-y)\dif s\preceq\varrho^0_{\alpha-1}(t,x-y).
$$
For $R_3(t,x,y)$, by  (\ref{ER6}), (\ref{eq3}) and (\ref{EU7}), we have
\begin{align*}
R_3(t,x,y)&\preceq\int^{\frac{t}{2}}_0\!\!\!\int_{\mR^d}\varrho^0_{\alpha-1}(t-s,x-z)(\varrho^0_\beta+\varrho^\beta_0)(s,z-y)\dif z\dif s\preceq\varrho^0_{\alpha-1}(t,x-y).
\end{align*}
Combining the above calculations and by (\ref{ER65}), we obtain (\ref{GH2}).
\end{proof}

Below, we write
$$
P^\kappa_t
f(x):=\int_{\mR^d}p^\kappa_{\alpha} (t,x,y)f(y)\dif y.
$$
\bl
For any bounded and H\"older continuous function $f$, we have
\begin{align}
\sL^\kappa_\alpha\left(\int^t_0 P^\kappa_sf(\cdot)\dif s\right)(x)=\int^t_0\sL^\kappa_\alpha P^\kappa_sf(x)\dif s, \ \ x\in\mR^d.\label{UY1}
\end{align}
\el
\begin{proof}
By definition of $\sL^\kappa_\alpha$ and Fubini's theorem, we have
\begin{align*}
\sL^\kappa_\alpha\left(\int^t_0 P^\kappa_sf\dif s\right)(x)
=\lim_{\eps\downarrow 0}\int_{|w|>\eps}\left(\int^t_0\delta_{ P^\kappa_sf}(x;w)\dif s\right)\kappa(x,w)|w|^{-d-\alpha}\dif w
=\lim_{\eps\downarrow 0}\int^t_0I_\eps(s,x)\dif s,
\end{align*}
where
$$
I_\eps(s,x):=\int_{|w|>\eps}\delta_{ P^\kappa_sf}(x;w)\kappa(x,w)|w|^{-d-\alpha}\dif w.
$$
Using the same argument as in proving (\ref{ETR1}), one can prove that for some $p>1$,
$$
\sup_{\eps\in(0,1)}\int^t_0|I_\eps(s,x)|^p\dif s<+\infty.
$$
Hence, we can interchange the limit and integral, and obtain (\ref{UY1}).
\end{proof}

\bl\label{Le43}
For any $p\in[1,\infty)$ and $f\in L^p(\mR^d)$, $(0,1)\ni t\mapsto\sL^\kappa_\alpha  P^\kappa_tf\in L^p(\mR^d)$ is continuous.
In the case of $p=\infty$, i.e., if $f$ is a bounded measurable function on $\mR^d$, then for each $x\in\mR^d$, $t\mapsto\sL^\kappa_\alpha  P^\kappa_tf(x)$
is a continuous function on $(0,1)$. Moreover, for any $p\in[1,\infty]$, there exists a constant
$C=C(p,d,\alpha,\beta, \kappa_0,\kappa_1, \kappa_2, p)>0$
such that for all $f\in L^p(\mR^d)$
and $t>0$,
\begin{align}
\|\sL^\kappa_\alpha  P^\kappa_tf\|_p\leq C t^{-1}\|f\|_p.\label{ERT1}
\end{align}
\el
\begin{proof}
For any $p\in[1,\infty]$, by Lemma \ref{Le01} and Young's inequality, we have
$$
\|\sL^\kappa_\alpha P^\kappa_tf\|_p\preceq\left(\int_{\mR^d}\left|\int_{\mR^d}\varrho^0_0(t,x-y)|f(y)|\dif y\right|^p\dif x\right)^{1/p}\leq\|\varrho^0_0(t)\|_1\|f\|_p\stackrel{(\ref{ES4})}{\preceq} t^{-1}\|f\|_p.
$$
Thus, we obtain (\ref{ERT1}).

On the other hand, for any $\eps\in(0,1)$, by Lemma \ref{Le01}, we have
$$
\sup_{t\in(\eps,1)}|\sL^\kappa_\alpha p^\kappa_{\alpha} (t,x,y)|\preceq\sup_{t\in(\eps,1)}\varrho^0_0(t,x-y)\preceq \varrho^0_0(\eps,x-y).
$$
Since for fixed $x\not=y\in\mR^d$, the mapping $t\mapsto\sL^\kappa_\alpha p^\kappa_{\alpha} (t,x,y)$ is continuous by (\ref{EP22}), the desired continuity of
$t\mapsto\sL^\kappa_\alpha  P^\kappa_{t}f(x)$ follows by the dominated convergence theorem.
\end{proof}

\subsection{Proof of Theorem \ref{Main}}
After the above preparation, we are now in a position to give the proof of Theorem \ref{Main}.
First of all, using Lemmas \ref{Le44} and \ref{Le6}, one sees that the calculations in (\ref{ER91}) make sense, and thus, we obtain (\ref{eq15}).

(i) Notice that by (\ref{ER66}), (\ref{eq3})  and (\ref{eq30}),
\begin{align}
\int^t_0\!\!\!\int_{\mR^d}p_z(t-s,x-z)|q(s,z,y)|\dif z\dif s
&\preceq\int^t_0\!\!\!\int_{\mR^d}\varrho^0_\alpha(t-s,x-z)(\varrho^0_\beta+\varrho^\beta_0)(s,z-y)\dif z\dif s\no\\
&\preceq (\varrho^0_{\alpha+\beta}+\varrho^\beta_\alpha)(t,x-y),\label{eq7}
\end{align}
which in turn gives estimate (\ref{eq16})
 by equation (\ref{ER65}) and (\ref{ER66}),
 where the constant $c_1$ can be chosen to depend only on
 $(d, \alpha, \beta, \kappa_0, \kappa_1, \kappa_2)$.

(ii) Estimate (\ref{Ho}) follows by Lemma \ref{Le00}.

(iii) Estimate (\ref{eq17}) follows by (\ref{GH1}). The continuity of $t\mapsto \sL^\kappa_\alpha p^\kappa_\alpha(t,\cdot,y)(x)$ follows by (\ref{EP22}).

(iv) Let $f$ be a bounded and uniformly continuous function. For any $\eps>0$, there exists a $\delta>0$ such that for all $|x-y|\leq\delta$,
$$
|f(x)-f(y)|\leq\eps.
$$
By (\ref{E32}) and (\ref{ER66}), we have
\begin{align*}
&\lim_{t\downarrow 0}\sup_{x\in\mR^d}\left|\int_{\mR^d}p_y(t,x-y)f(y)\dif y-f(x)\right|\no\\
&\quad\preceq \lim_{t\downarrow 0}\sup_{x\in\mR^d}\int_{\mR^d}\varrho^0_\alpha(t,x-y)\cdot |f(y)-f(x)|\dif y\no\\
&\quad\preceq \eps+2\|f\|_\infty\lim_{t\downarrow 0}\sup_{x\in\mR^d}\int_{|x-y|>\delta}\varrho^0_\alpha(t,x-y)\dif y\leq\eps,
\end{align*}
which implies that
\begin{align*}
\lim_{t\downarrow 0}\sup_{x\in\mR^d}\left|\int_{\mR^d}p_y(t,x-y)f(y)\dif y-f(x)\right|=0.
\end{align*}
Moreover, by (\ref{eq7}), we also have
\begin{align*}
&\left|\int_{\mR^d}\!\int_0^t\!\!\!\int_{\mR^d}p_z(t-s,x-z)q(s,z,y)f(y)\dif z\dif s\dif y\right|\\
&\quad\preceq \int_{\mR^d}(\varrho^0_{\alpha+\beta}+\rho^\beta_\alpha)(t,x-y)\dif y
\stackrel{(\ref{ES4})}{\preceq} t^{\frac{\beta}{\alpha}}\rightarrow 0,\ t\downarrow 0,
\end{align*}
Thus, (\ref{eq14}) is proven by equation (\ref{ER65}).

\vspace{5mm}

We now show that  kernels that satisfy \eqref{eq15}--\eqref{eq14} is unique.
For this, let $\widetilde p^\kappa_\alpha (t, x, y)$ be  any kernel that satisfies \eqref{eq15}--\eqref{eq14} and, for
 $f\in C^\infty_c(\mR^d)$, define $\wt u_f(t,x):=\int_{\mR^d}
 \wt p^\kappa_\alpha (t, x, y) f(y) \dif y$.
 First of all, by (iv), one sees that
$$
\wt u_f\in C_b([0,1]\times\mR^d),\quad
 \lim_{t\downarrow 0}\sup_{x\in\mR^d}| \wt u_f(t,x)-f(x)|=0.
$$
Secondly, by
\eqref{Ho}
we have for any $\gamma\in(0,\alpha\wedge 1)$,
\begin{align*}
|\wt u_f(t,x)- \wt u_f(t,x')|&\leq\|f\|_\infty\int_{\mR^d}
|\wt p^\kappa_{\alpha} (t,x,y)-\wt p^\kappa_{\alpha} (t,x',y)|\dif y\\
&\preceq\|f\|_\infty|x-x'|^\gamma\int_{\mR^d}
\left(\varrho^0_{\alpha-\gamma}(t,x-y)+\varrho^0_{\alpha-\gamma}(t,x'-y)\right)
\dif y\\
&\stackrel{(\ref{ES4})}{\preceq}\|f\|_\infty|x-x'|^\gamma t^{-\frac{\gamma}{\alpha}}.
\end{align*}
The same holds for $u_f (t, x):=  \int_{\mR^d}
   p^\kappa_\alpha (t, x, y) f(y) \dif y$.
Thus in view of \eqref{eq15} and (iii),
$w(t, x):= u_f(t, x)-\wt u_f (t, x)$ satisfies
all the conditions of Theorem \ref{Th1} with $w(0, x)=0$ for every $x\in \mR^d$.
Applying Theorem \ref{Th1} to both $w$ and $-w$ yields
$w(t, x)=0$ for every $t>0$ and $x\in \mR^d$.
Consequently, we have $\wt p^\kappa_\alpha (t, x,y)=p^\kappa_\alpha (t, x,y)$.

\vspace{5mm}

(1) has already been proved in the above.

(2) Applying the maximum principle Theorem \ref{Th1} to $u_f$ with
$f\in C^\infty_c (\mR^d)$ and $f\leq 0$ implies that
$p^\kappa_\alpha (t, x, y)\geq 0$. Moreover, since constant function
$u(t, x)=1$ solves the equation $\partial_t u(t, x)=\sL^\kappa_\alpha u(t, x)$
with initial value 1, we have \eqref{EK1}.

(3) This follows from the uniqueness of the solution to $\partial_t u(t, x)=\sL^\kappa_\alpha u(t, x)$, implied by Theorem \ref{Th1}.

(4) will be proven in the next subsection.

(5) If $\alpha\in[1,2)$, then estimate (\ref{eq170}) follows by (\ref{GH2}).

(6) For $f\in C^2_b(\mR^d)$, define
$$
u(t,x):=f(x)+\int^t_0 P^\kappa_s\sL^\kappa_\alpha f(x)\dif s.
$$
By (\ref{UY1}) we have
\begin{align*}
\sL^\kappa_\alpha u(t,x)&=\sL^\kappa_\alpha f(x)+\int^t_0\sL^\kappa_\alpha P^\kappa_s\sL^\kappa_\alpha f(x)\dif s\\
&=\sL^\kappa_\alpha f(x)+\int^t_0\p_s( P^\kappa_s\sL^\kappa_\alpha f)(x)\dif s\\
&= P^\kappa_t\sL^\kappa_\alpha f(x)=\p_t u(t,x).
\end{align*}
Moreover, it is easy to see that (\ref{ER90}), (\ref{ER93}) and (\ref{ER92}) are satisfied for $u$. Thus, by Theorem \ref{Th1} we obtain
\begin{equation}\label{e:4.23}
 P^\kappa_t f(x)=u(t,x)=f(x)+\int^t_0 P^\kappa_s\sL^\kappa_\alpha f(x)\dif s,
\end{equation}
which in turn implies that
$$
\lim_{t\downarrow 0}\frac{1}{t}( P^\kappa_t f(x)-f(x))=\lim_{t\downarrow 0}\frac{1}{t}\int^t_0 P^\kappa_s\sL^\kappa_\alpha f(x)\dif s\stackrel{(\ref{eq14})}{=}\sL^\kappa_\alpha f(x)
$$
and the convergence is uniform.

(7) Fix $p\in[1,\infty)$. By (iv), (2) and (\ref{ERT1}), it is easy to see that
$( P^\kappa_t)_{t\geq 0}$ is a $C_0$-semigroup in $L^p(\mR^d)$. On the other hand, for any $f\in L^p(\mR^d)$,
by equation (\ref{eq15}) and Lemma \ref{Le43}, one sees that $ P^\kappa_t f$ is differentiable in $L^p(\mR^d)$ for any $t>0$, i.e.,
\begin{align*}
&\lim_{\eps\to 0}\frac{\| P^\kappa_{t+\eps} f- P^\kappa_t f-\eps\sL^\kappa_\alpha P^\kappa_t f\|_p}{\eps}\\
&\quad\leq\lim_{\eps\to 0}\frac{1}{\eps}\int^{t+\eps}_t\|\sL^\kappa_\alpha P^\kappa_{t+s} f-\sL^\kappa_\alpha P^\kappa_t f\|_p\dif s=0.
\end{align*}
The analyticity of $C_0$-semigroup $( P^\kappa_t)_{t\geq 0}$
follows by (\ref{ERT1}) and \cite[p.61 Theorem 5.2 (d)]{Pa}.

\subsection{Proof of lower bound estimate of $p^\kappa_\alpha(t,x,y)$}

From the previous subsection, one sees that $( P^\kappa_t)_{t\geq 0}$ is a Feller semigroup.
Hence, it determines a Feller process $(\Omega,\sF, (\mP_x)_{x\in\mR^d}, (X_t)_{t\geq 0})$.
For any $f\in C^2_b(\mR^d)$,
it follows from \eqref{e:4.23} and the Markov property of $X$
that under $\mP_x$, with respect to the filtration $\sF_t:=\sigma\{X_s, s\leq t\}$,
\begin{align}
M^f_t:=f(X_t)-f(X_0)-\int^t_0\sL^\kappa_\alpha f(X_s)\dif s\ \mbox{ is a martingale}.\label{ERY1}
\end{align}
In other words, $\mP_x$ solves the martingale problem for $(\sL^\kappa_\alpha,
C^2_b (\mR^d))$. Thus $\mP_x$ in particular solves the martingale problem
for $(\sL^\kappa_\alpha, C^\infty_c (\mR^d))$.

We now derive a L\'evy system of $X$ by following an approach from \cite{Ch-Ki-So}.
By (\ref{ERY1}), one can derive that $X_t=(X^{1}_t, \dots, X^{d}_t) $ is a semi-martingale. By It\^o's formula,
we have that, for any $f\in C^\infty_c(\mR^d)$,
\begin{equation}\label{e:ito}
f(X_t)-f(X_0)=\sum^d_{i=1}\int^t_0{\partial}_if(X_{s-})\dif X^{i}_s +\sum_{s\leq t}\eta_s(f) +\frac12 A_t(f),
\end{equation}
where
\begin{equation}\label{e:ito2}
\eta_s(f)=f(X_s)-f(X_{s-})-\sum^d_{i=1}{\partial}_if(X_{s-})(X^{i}_s-X^{i}_{s-})
\end{equation}
and
\begin{equation}\label{e:ito3}
A_t(f)=\sum^d_{i, j=1}\int^t_0{\partial}_i{\partial}_jf(X_{s-})\dif\langle (X^{i})^c, (X^{j})^c\rangle_s.
\end{equation}

Now suppose that $A$ and $B$ are two bounded closed subsets of $\mR^d$
having a positive distance from each other. Let $f\in C^\infty_c(\mR^d)$ with
$f=0$ on $A$ and $f=1$ on $B$.  Clearly
$N^f_t:=\int^t_0{\bf 1}_A(X_{s-})\dif M^f_s$ is a martingale.
Define
\begin{equation}
 J(x, y)= k(x, y-x)/|y-x|^{d+\alpha},
\end{equation}\label{e:4.28}
so $\sL^\kappa_\alpha$ can be rewritten as
\begin{equation}\label{e:L2}
\sL^\kappa_\alpha f(x)= \lim_{\eps\to 0} \int_{\{|y-x|>\eps\}}
(f(y)-f(x)) J(x, y) \dif y.
\end{equation}
We get by  \eqref{ERY1}--\eqref{e:ito3} and  \eqref{e:L2},
\begin{align*}
N^f_t&=\sum_{s\leq t}{\bf 1}_A(X_{s-})(f(X_s)-f(X_{s-})) -\int^t_0{\bf 1}_A(X_s)\sL^\kappa_\alpha f(X_s)\dif s\\
&=\sum_{s\leq t}{\bf 1}_A(X_{s-})f(X_s)-\int^t_0{\bf 1}_A(X_s)\int_{\mR^d}f(y)J(X_s,y)\dif y\dif s.
\end{align*}
By taking a sequence of functions $f_n\in C^\infty_c(\mR^d)$ with
$f_n=0$ on $A$,  $f_n=1$ on $B$ and $f_n\downarrow {\bf 1}_B$, we get
that, for any $x\in \mR^d$,
$$
\sum_{s\leq t}{\bf 1}_A(X_{s-}){\bf 1}_B(X_s) -\int^t_0{\bf1}_A(X_s)\int_BJ(X_s, y)\dif y\dif s
$$
is a martingale with respect to $\mP_x$. Thus,
$$
\mE_x\left[ \sum_{s\leq t}{\bf 1}_A(X_{s-}){\bf 1}_B(X_s)\right]=
\mE_x\left[\int^t_0\int_{\mR^d} {\bf 1}_A(X_s){\bf 1}_B(y)J(X_s, y)\dif y\dif s\right].
$$
Using this and a routine measure theoretic argument, we get
$$
\mE_x\left[ \sum_{s\leq t}f(X_{s-}, X_s) \right]
=\mE_x\left[\int^t_0\int_{\mR^d}f(X_s, y)J(X_s, y)\dif y\dif s\right]
$$
for any non-negative measurable function $f$ on $\mR^d\times \mR^d$
vanishing on $\{(x, y)\in \mR^d\times \mR^d: x=y\}$. Finally
following the same arguments as in \cite[Lemma 4.7]{Ch-Ku} and
\cite[Appendix A]{Ch-Ku1}, we get

\bt\label{T:l2}
$X$ has  a L\'evy system $(J, t)$ as $X$, that is, for any
$x\in \mR^d$ and any non-negative measurable function $f$ on $\mR_+
\times \mR^d\times \mR^d$ vanishing on $\{(s, x, y)\in \mR_+ \times
\mR^d\times \mR^d: x=y\}$ and $(\sF_t)$-stopping time $T$,
\begin{equation}\label{e:ls4xb}
\mE_x \left[\sum_{s\leq T} f(s,X_{s-}, X_s) \right]= \mE_x \left[
\int_0^T \left( \int_{\mR^d} f(s,X_s, y) J(X_s, y)\dif y \right) \dif s
\right].
 \end{equation}
\et

For a set $K\subset\mR^d$, denote
$$
\sigma_K:=\inf\{t\geq 0: X_t\in K \},\ \ \tau_K:=\inf\{t\geq 0: X_t\notin K\}.
$$
Let $B(x,r)$ be the ball with radius $r$ and center $x$. We need the following lemma (see \cite{Ba-Le, Ch-Ku}).

\bl\label{Le56}
For each $\gamma\in(0,1)$, there exists $A_0>0$ such that for every $A>A_0$ and $r\in(0,1)$,
\begin{align}
\mP_x(\tau_{B(x,Ar)} \leq r^\alpha)\leq \gamma.\label{ERY3}
\end{align}
\el

\begin{proof}
Without loss of generality, we assume that $x=0$.
Given $f\in C^2_b(\mR^d)$ with $f(0)=0$ and $f(x)=1$ for $|x|\geq 1$,  we set
$$
f_r(x):=f(x/r),\ \ r>0.
$$
By the definition of $f_r$, we have
\begin{align}
 \mP_0(\tau_{B(0,Ar)}\leq r^\alpha)
\leq \mE_0 \left[ f_{Ar} ( X_{\tau_{B(0,Ar)}\wedge r^{\alpha}}) \right]
\stackrel{(\ref{ERY1})}{=}\mE_0\left(\int^{\tau_{B(0,Ar)}\wedge r^{\alpha}}_0\sL^\kappa_\alpha f_{Ar}(X_s)\dif s\right).\label{ERY2}
\end{align}
On the other hand, by the definition of $\sL^\kappa_\alpha $, we have for $\lambda>0$,
\begin{align*}
|\sL^\kappa_\alpha f_{Ar}(x)|&=\frac{1}{2}\left|\int_{\mR^d}(f_{Ar}(x+z)+f_{Ar}(x-z)-2f_{Ar}(x))\kappa(x,z)|z|^{-d-\alpha}\dif z\right|\\
&\leq\frac{\kappa_1\|\nabla^2 f_{Ar}\|_\infty}{2}\int_{|z|\leq\lambda r}|z|^{2-d-\alpha}\dif z+2\kappa_1\|f_{Ar}\|_\infty\int_{|z|\geq\lambda r}|z|^{-d-\alpha}\dif z\\
&=\kappa_1\frac{\|\nabla^2 f\|_\infty}{(Ar)^2}\frac{(\lambda r)^{2-\alpha}}{2(2-\alpha)}s_1+2\kappa_1\|f\|_\infty\frac{(\lambda r)^{-\alpha}}{\alpha}s_1\\
&=\kappa_1s_1\left(\frac{\|\nabla^2 f\|_\infty}{A^2}\frac{\lambda^{2-\alpha}}{2(2-\alpha)}+2\|f\|_\infty
\frac{\lambda^{-\alpha}}{\alpha}\right) r^{-\alpha},
\end{align*}
where $s_1$ is the sphere area of the unit ball. Substituting this into (\ref{ERY2}),
we get
$$
\mP_0(\tau_{B(0,Ar)}\leq r^\alpha) \leq
\kappa_1s_1\left(\frac{\|\nabla^2 f\|_\infty}{A^2}\frac{\lambda^{2(2-\alpha)}}{2-\alpha}+2\|f\|_\infty
\frac{\lambda^{-\alpha}}{\alpha}\right) .
$$
Choosing
first $\lambda$ large enough and then $A$ large enough
yield the desired estimate.
\end{proof}
Now we can give
\begin{proof}[Proof of lower bound of $p^\kappa_\alpha(t,x,y)$]
By Lemma \ref{Le56}, there is a constant $\lambda\in (0,\tfrac{1}{2})$ such that for all $t\in(0,1)$,
\begin{align}
\mP_x(\tau_{B(x,t^{1/\alpha}/2)}\leq \lambda t)\leq\tfrac{1}{2}.\label{ERY7}
\end{align}
By (\ref{ER65}), (\ref{ER66}) and (\ref{eq7}), there is a time $t_0\in(0,1)$ such that
$$
p^\kappa_\alpha(t,x,y)\succeq t^{-d/\alpha}\mbox{ for all $t\in(0,t_0)$ and $|x-y|\leq 3t^{1/\alpha}$}.
$$
By C-K equation \eqref{eq21} and iterating $[1/t_0]+1$ times, we conclude that
\begin{align}
p^\kappa_\alpha(t,x,y)\succeq t^{-d/\alpha}\mbox{ for all $t\in(0,1)$ and $|x-y|\leq 3t^{1/\alpha}$}.\label{ERY9}
\end{align}
Below, we assume
\begin{align}
|x-y|>3t^{1/\alpha}.\label{ERY8}
\end{align}
For the given $\lambda$ in (\ref{ERY7}), by the strong Markov property, we have
\begin{align}
\mP_x\left(X_{\lambda t}\in B(y,t^{1/\alpha})\right)
&\geq \mP_x\left(\sigma:=\sigma_{B(y,t^{1/\alpha}/2)}
\leq \lambda t; \sup_{s\in[\sigma,\sigma+\lambda t]}|X_s-X_\sigma|< t^{1/\alpha}/2\right)\no\\
&=\mE_x\left(\mP_z\left(\sup_{s\in[0,\lambda t]}|X_s-z| < t^{1/\alpha}/2\right)\Big|_{z=X_\sigma}; \sigma_{B(y,t^{1/\alpha}/2)} \leq \lambda t\right)\no\\
&\geq \inf_{z\in B(y,t^{1/\alpha}/2)}
\mP_z\left(\tau_{B(z,t^{1/\alpha}/2)}> \lambda t\right)
\mP_x\left(\sigma_{B(y,t^{1/\alpha}/2)}\leq \lambda t\right)\no\\
&\stackrel{(\ref{ERY7})}{\geq}\tfrac{1}{2}
\mP_x\left(\sigma_{B(y,t^{1/\alpha}/2)} \leq  \lambda t\right)\no\\
&\geq \tfrac{1}{2}\mP_x\left(X_{\lambda t\wedge\tau_{B(x,t^{1/\alpha})}}\in B(y,t^{1/\alpha}/2)\right).\label{ERY11}
\end{align}
Noticing that
$$
X_s\notin B(y,t^{1/\alpha}/2)\subset B(x,t^{1/\alpha})^c,\ \ s<\lambda t\wedge\tau_{B(x,t^{1/\alpha})},
$$
we have
$$
{\bf 1}_{X_{\lambda t\wedge\tau_{B(x,t^{1/\alpha})}}\in B(y,t^{1/\alpha}/2)}=\sum_{s\leq \lambda t\wedge\tau_{B(x,t^{1/\alpha})}}{\bf 1}_{X_s\in B(y,t^{1/\alpha}/2)}.
$$
Thus, by (\ref{e:ls4xb}) we have
\begin{align}
\mP_x\left[ X_{\lambda t\wedge\tau_{B(x,t^{1/\alpha})}}\in B(y,t^{1/\alpha}/2)\right]
&=\mE_x\left[\int^{\lambda t\wedge\tau_{B(x,t^{1/\alpha})}}_0\!\!\int_{B(y,t^{1/\alpha}/2)} J(X_s,u)\dif u\dif s\right]\no\\
&\geq\mE_x\left[\int^{\lambda t\wedge\tau_{B(x,t^{1/\alpha})}}_0\!\!\int_{B(y,t^{1/\alpha}/2)}\frac{\kappa_0}{|X_s-u|^{d+\alpha}}\dif u\dif s\right]\no\\
&\geq\mE_x\left[\lambda t\wedge\tau_{B(x,t^{1/\alpha})}\right]\int_{B(y,t^{1/\alpha}/2)}
 \frac{\kappa_0}{(|x-y|+3t^{1/\alpha}/2)^{d+\alpha}}\dif u\no\\
&\stackrel{(\ref{ERY8})}{\geq} \lambda t\, \mP_x\left(\tau_{B(x,t^{1/\alpha})}\geq\lambda t\right) \left(\int_{B(y,t^{1/\alpha}/2)}\dif u\right) \frac{\kappa_0(2/3)^{d+\alpha}}{|x-y|^{d+\alpha}}\no\\
&\geq \Big(\lambda \kappa_0(2/3)^{d+\alpha}2^{-d-1}s_1\Big)\frac{t^{1+d/\alpha}}{|x-y|^{d+\alpha}},\label{ERY10}
\end{align}
where $s_1$ is the sphere area of the unit ball.

Now, by Chapman-Kolmogorov's equation again, we have
\begin{align*}
p^\kappa_\alpha(t,x,y)&\geq \int_{B(y,t^{1/\alpha})}p^\kappa_\alpha(\lambda t,x,z)p^\kappa_\alpha((1-\lambda) t, z,y)\dif z\\
&\geq \inf_{z\in B(y,t^{1/\alpha})}p^\kappa_\alpha((1-\lambda) t, z,y)\int_{B(y,t^{1/\alpha})}p^\kappa_\alpha(\lambda t,x,z)\dif z\\
&\stackrel{(\ref{ERY9})}{\succeq} t^{-d/\alpha}\mP_x\left(X_{\lambda t}\in B(y,t^{1/\alpha})\right)\stackrel{(\ref{ERY11}),(\ref{ERY10})}{\succeq} t|x-y|^{-d-\alpha}.
\end{align*}
which, combining with (\ref{ERY9}), gives the lower bound estimate of $p^\kappa_\alpha(t,x,y)$.
\end{proof}

\subsection{Proof of Corollary \ref{C:1.4}}
Since $\lambda_0 I_{d\times d} \leq A(x)) \leq \lambda_1 I_{d\times d}$
and $|a_{ij}(x)-a_{ij}(y)|\leq \lambda_2 |x-y|^\beta$ for each
$1\leq i,j\leq d$, the function
$\kappa (x, z)$ defined by \eqref{e:1.21} satisfies the conditions
\eqref{Con1}-\eqref{Con2} with $\kappa_i$, $i=0, 1, 2$, depend only
on $d, \alpha, \beta, \lambda_0, \lambda_1$ and $\lambda_2$.
Thus by Theorem \ref{Main}, there is a jointly continuous heat kernel
$p(t, x, y)$ for the non-local operator $\sL=\sL^\kappa_\alpha$
of \eqref{e:1.20} corresponding to this $\kappa (x, z)$.
Let
$\{\wt X, \mP_x, x\in \mR^d\}$
be the Feller process having
$p(t, x , y)$ as its transition density function.
 As we observed in the beginning of subsection \S 4.4,
 $\mP_x$ solves the martingale problem for $(\sL, C^2_b(\mR^d))$.
 On the other hand,
 it is shown in \S 7 of \cite{BC} (see Theorem 7.1 and its proof
 as well as Theorems 4.1 and  6.3 there)
 that the law of
 the unique weak solution $X$ to SDE \eqref{e:1.18} is the unique
 solution to the martingale problem for $(\sL, C^2_b(\mR^d))$.
 Hence $\wt X$ and $X$ have the same distribution.
 Therefore $p(t, x, y)$ is the transition density function
 of $X$. The conclusion of the corollary now follows from Theorem \ref{Main}.
 \qed

\bigskip

{\bf Acknowledgements:}
The authors are grateful to
Renming Song  and   Longjie Xie
for their useful conversations
over an earlier version of this paper.


\begin{thebibliography}{10}

\bibitem{BC} Bass R.F. and Chen Z.-Q.:
Systems of equations driven by stable processes. {\it Probab. Theory Relat. Fields \bf 134} (2006), 175--214.


\bibitem{Ba-Le}Bass R.F. and Levin D.A.: Harnack inequalities for jump processes. {\it Potential Anal. \bf 17} (2002), 375-388.

\bibitem{BL2}Bass R.F. and Levin D.A.: Transition probabilities for symmetric
jump processes. {\it Trans. Amer. Math. Soc. \bf 354} (2002), 2933-2953.

\bibitem{BR} Bass, R.F. and Ren, H.: Meyers inequality and strong stability for stable-like operators. {\it J. Funct. Anal. \bf 265} (2013), 28--48.


\bibitem{Bl-Ge}Blumenthal R.M. and Getoor R.K.: Some theorems on stable processes. {\it Trans. Amer. Math. Soc. \bf 95} (1960), 263-273.

\bibitem{Bo-Ja0}Bogdan K., Jakubowski T. : Estimates of heat kernel of fractional Laplacian perturbed by gradient operator.
{\it Commun. Math. Phys., \bf 271} (2007), 179-198.

\bibitem{Bo-Ja1}Bogdan K., Jakubowski T. : Estimates of Green function for the fractional Laplacian perturbed by gradient. {\it Potential Anal. \bf 36}
    (2012), 455-481.

\bibitem{CaS1} Caffarelli L. and Silvestre L.:
The Evans-Krylov theorem for nonlocal fully nonlinear equations.
{\it Ann. Math. \bf 174} (2011), 1163-1187.


\bibitem{Ch}Chen Z.-Q.: Symmetric jump processes and their heat kernel esitmates. {\it Science in China Series A: Mathematics  \bf 52} (2009),  1423-1445.

\bibitem{CH} Chen Z.-Q. and Hu E.: Heat kernels of $\Delta+\Delta^{\alpha/2}$
perturbed by gradient operators.
In preparation.


\bibitem{Ch-Ki-So0}Chen Z.-Q., Kim P. and Song R.: Heat kernel estimates for Dirichlet fractional Laplacian. {\it J. European Math. Soc, \bf 12} (2010), 1307-1329.

\bibitem{Ch-Ki-So}Chen Z.-Q., Kim P. and Song R.: Dirichlet heat kernel estimates for fractional Laplacian with gradient perturbation.
{\it Ann. Probab. \bf 40} (2012), 2483-2538.

\bibitem{Ch-Ki-So1}Chen Z.-Q., Kim P. and Song R.: Stability of Dirichlet heat kernel estimates for non-local operators
under Feynmman-Kac perturbation. To appear in {\it Trans. Amer. Math.Soc.}

\bibitem{Ch-Ku}Chen Z.-Q. and Kumagai T.: Heat kernel estimates for stable-like processes on $d$-sets. {\it Stochastic Processes and their Applications
    \bf 108} (2003), 27-62.

\bibitem{Ch-Ku1}Chen Z.-Q. and   Kumagai T.: Heat kernel estimates for jump processes of mixed types on metric measure spaces. {\it Probab. Theory Relat. Fields  \bf 140} (2008), 277--317.

\bibitem{CW} Chen Z.-Q. and Wang J.: Perturbation by non-local operators.
Preprint, 2012.

\bibitem{CWa} Chen Z.-Q. and Wang L.: Heat kernel estimates for relativistic stable processes with singular drifts. In preparation, 2013.

\bibitem{CZh} Chen Z.-Q. and Zhang X.: Uniqueness of stable-like processes.
In preparation, 2013.

\bibitem{CZ} Chen Z.-Q. and Zhao Z.: Diffusion processes and second order elliptic operators with singular coefficients for lower order
   terms. {\it Math. Ann. \bf 302} (1995), 323-357.

\bibitem{De}Doetsch, G.: {\it Introduction to the Theory and Application of the Laplace Transformation}. New-York, Springer-Verlag, 1974.

\bibitem{Ei-Iv-Ko}Eidelman S. D.,  Ivasyshen S. D. and  Kochubei A. N.:
{\it Analytic Methods in the Theory of Differential and Pseudo-differential Equations of Parabolic Type}. Birkhauser, Basel, 2004.

\bibitem{Fr0}Friedman, A.: {\it Partial Differential Equations}.
Holt, Rinehart and Winston, Inc., New York-Montreal, Que-London, 1969. vi+262 pp.

\bibitem{Fr}Friedman A.: {\it Partial Differential Equations of Parabolic Type}. Prentice-Hall, Englewood Cliffs, N.J., 1975.

\bibitem{He}Henry D.: {\it Geometric Theory of Semilinear Parabolic Equations}. Lect. Notes in Math., Vol. 840, 1981.

\bibitem{Ja}Jakubowski, T.: Fractional Laplacian with singular drift.
{\it Studia Math. \bf 207} (2011), 257-273.

\bibitem{Ja-Sz1}Jakubowski T. and Szczypkowski K.: Estimates of gradient perturbation series. {\it J. Math. Anal. Appl. \bf 389} (2012), 452-460.

\bibitem{Ja-Sz2}Jakubowski T. and Szczypkowski K.: Time-dependent gradient perturbations of fractional Laplacian. {\it J. Evol. Equ. \bf 10} (2010), 319-339.

\bibitem{Ko}Kochubei A. N.: Parabolic pseudodifferential equations, hypersingular integrals and Markov processes. {\it Math. USSR, Izv. \bf 33} (1989), 233-259;
translation from {\it Izv. Akad. Nauk SSSR, Ser. Mat. \bf 52} (1988), 909-934.

\bibitem{Kol}
Kolokoltsov V.: Symmetric stable laws and stable-like jump-diffusions.
  {\it Proc. Lond. Math. Soc. \bf 80} (2000), 725-768.

\bibitem{La-So-Ur}Lady$\breve{z}$enskaja O.A., Solonnikov V.A., Ural'ceva N.N.: {\it Linear and Quasi-linear Equations of Parabolic Type}.
     Translated from the Russian by S.Smith. American Mathematical Society, 1968.

\bibitem{Ma-Mi}Maekawa, Y. and Miura, H.: Upper bounds for fundamental solutions to non-local diffusion equations with divergence free drift.
    {\it J. Funct. Anal. \bf 264} (2013), 2245-2268.


\bibitem{Pa}Pazy A.: {\it Semigroups of Linear Operators and Applications to Partial Differential Equations}. Springer-Verlag, New York, Berlin, 1983.

\bibitem{Wa-Zh} Wang F.Y. and Zhang X.: Heat kernel for fractional diffusion operators with perturbations. http://arxiv.org/abs/1204.4956.
To appear in {\it Forum Mathematicum}.

\bibitem{Xi-Zh}Xie L. and Zhang X.: Heat kernel estimates for critical fractional diffusion operator. http://arxiv.org/abs/1210.7063.

\bibitem{Zh1}Zhang X.: $L^p$-maximal regularity of nonlocal parabolic equation and applications.
{\it Annales de l'Institut Henri Poincare Analyse non lineaire. \bf 30} (2013), 573-614.

\bibitem{Zh2}Zhang X.: $L^p$-solvability of nonlocal parabolic equations with spatial dependent and non-smooth kernels. In {\it Emerging Topics on Differential Equations and their Applications}. Edited by Hua Chen, Yiming Long and Yasumasa Nishiura, World Scientific, 247-262 (2013), http://arxiv.org/abs/1206.2709.

\end{thebibliography}
\end{document}